\documentclass[journal,10pt,twoside]{IEEEtran}  

\IEEEoverridecommandlockouts                              
\makeatletter
\newcommand{\removelatexerror}{\let\@latex@error\@gobble}
\makeatother

\usepackage[english]{babel}
\usepackage[usenames,dvipsnames,table]{xcolor}
\usepackage{amsmath}
\usepackage{amsthm}
\usepackage{amssymb}
\usepackage{graphicx}
\usepackage[footnotesize]{caption}
\usepackage[font=footnotesize]{subcaption}
\usepackage{url}
\usepackage{float}
\usepackage{tikz}

\usetikzlibrary{shapes,arrows,positioning}

\usepackage{pgfplots} 		       
\usepackage[utf8]{inputenc}
\usepackage[hidelinks]{hyperref}
\usepackage{dsfont}   		       
\usepackage{siunitx}
\usepackage[mathscr]{euscript}
\usepackage{mathtools}
\usepackage{xargs} 		       
\usepackage[ruled,shortend,commentsnumbered,linesnumbered]{algorithm2e}

\usepackage{etoolbox}
\let\bbordermatrix\bordermatrix
\patchcmd{\bbordermatrix}{8.75}{4.75}{}{}
\patchcmd{\bbordermatrix}{\left(}{\left[}{}{}
    \patchcmd{\bbordermatrix}{\right)}{\right]}{}{}

\graphicspath{{./img/}{./}} 

\newcommand{\w}{\omega}
\newcommand{\e}{\epsilon}
\newcommand{\1}{\mathds 1  }

\newcommand{\vp}{\varphi}
\newcommand{\ul}{\underline}
\newcommand{\ol}{\overline}
\newcommand{\n}{\nabla}
\newcommand{\G}{\mathcal{G} }
\newcommand{\V}{\mathcal{V} }

\newcommand{\E}{\mathcal{E} }

\newcommand{\D}{\mathscr{D} }

\renewcommand{\l}{\lambda }

\newcommand{\ddt}[1][]{\frac{d #1}{dt}}
\newcommand{\myint}[1][\textstyle]{#1 \int}
\newcommand{\mysum}[1][\textstyle]{#1 \sum}

\newcommand{\gwin}{Gronwall's inequality}

\newcommand{\zetaobar}{\ol \zeta}
\newcommand{\xiobar}{\ol \xi}
\newcommand{\zetaubar}{\ul \zeta}
\newcommand{\xiubar}{\ul \xi}

\newcommand{\baalgo}{\textsc{Iterative Bid Update and Market Clearing Algorithm}\xspace}

\newcommand{\coeffc}{c}

\newcommand{\gain}{T}
\newcommand{\gainl}{\tau}

\newcommandx{\F}[2][1=1, 2=2]{\frac{#1}{#2}} 

\DeclareMathOperator{\col}{col}

\DeclareMathOperator{\diag}{diag}

\DeclareMathOperator{\blockdiag}{blockdiag}
\DeclareMathOperator*{\argmax}{arg\,max}

\newcommand{\real}{\mathbb{R}}
\newcommand{\realnonnegative}{{\mathbb{R}}_{\ge 0}}
\newcommand{\norm}[1]{\ensuremath{\| #1 \|}}
\newcommand{\norminf}[1]{\ensuremath{\| #1 \|_\infty}}
\newcommand{\until}[1]{[#1]}
\newcommand{\longthmtitle}[1]{\mbox{}\textup{\textsl{(#1):}}}

\allowdisplaybreaks

\newcommand{\st}{\operatorname{subject \text{$\, \,$} to}}
\newcommand{\minimize}{\operatorname{minimize}}
\newcommand{\oprocendsymbol}{\hbox{$\bullet$}}
\newcommand{\oprocend}{\relax\ifmmode\else\unskip\hfill\fi\oprocendsymbol}

\usepackage[normalem]{ulem} 
\definecolor{new}{rgb}{0.55,0,0.55}

\makeatletter
\newcommand{\thickhline}{%
  \noalign {\ifnum 0=`}\fi \hrule height 1pt
  \futurelet \reserved@a \@xhline
}
\newcolumntype{"}{@{\hskip\tabcolsep\vrule width 1pt\hskip\tabcolsep}}
\makeatother

\newtheorem{theorem}{Theorem}[section]
\newtheorem{proposition}[theorem]{Proposition}

\theoremstyle{remark}

\theoremstyle{definition}

\usepackage[prependcaption,colorinlistoftodos]{todonotes}

\title{\Large \bf Hybrid interconnection of iterative bidding and
  power network dynamics\\
  for frequency regulation and optimal dispatch
\thanks{A preliminary version of this work appeared
    as~\cite{ACC18-TS-AC-CDP-AS-JC} at the American Control
    Conference.}}

\author{Tjerk~Stegink,~Ashish~Cherukuri,~Claudio~De~Persis,~Arjan~van~der~Schaft,~and~Jorge
  Cort\'es
  \thanks{This work is supported by the
    NWO 
    \emph{Uncertainty Reduction in Smart
      Energy Systems (URSES)} program and the ARPA-e \emph{Network
      Optimized Distributed Energy Systems (NODES)} program.}
  \thanks{T. W. Stegink, C. De Persis and A. J. van der Schaft are
    with the Jan C. Willems Center for Systems and Control, University
    of Groningen, 9747 AG Groningen, the Netherlands.  {\tt\small
      \{t.w.stegink, c.de.persis,
      a.j.van.der.schaft\}@rug.nl}}
      \thanks{A. Cherukuri is with the Automatic Control Laboratory,
	ETH 
	Z\"{u}rich. {\tt\small cashish@control.ee.ethz.ch}}
   \thanks{J. Cort\'es is with the Department of
    Mechanical and Aerospace Engineering, University of California,
    San Diego.  {\tt\small cortes@ucsd.edu}}}

\begin{document}

\maketitle
\thispagestyle{empty}
\pagestyle{empty}

\begin{abstract}
  This paper considers a real-time electricity market involving an
  independent system operator (ISO) and a group of strategic
  generators. The ISO operates a market where generators bid prices at
  which there are willing to provide power. The ISO makes power
  generation assignments with the goal of solving the economic
  dispatch problem and regulating the network frequency.  We propose a
  multi-rate hybrid algorithm for bidding and market clearing that
  combines the discrete nature of iterative bidding with the
  continuous nature of the frequency evolution in the power network.
  We establish sufficient upper bounds on the inter-event times that
  guarantee that the proposed algorithm asymptotically converges to an
  equilibrium corresponding to an efficient Nash equilibrium and zero
  frequency deviation. Our technical analysis builds on the
  characterization of the robustness properties of the continuous-time
  version of the bidding update process interconnected with the power
  network dynamics via the identification of a novel LISS-Lyapunov
  function.  Simulations on the IEEE 14-bus system illustrate our
  results.
\end{abstract}

\section{Introduction}


The dispatch of power generation in the grid has been traditionally
done in a hierarchical fashion. Broadly speaking, cost efficiency is
ensured via market clearing at the upper layers and frequency
regulation is achieved via primary and secondary controllers at the
bottom layers. Research on improving the performance of these layers
has mostly developed independently from each other, motivated by their
separation in time-scales.  The increasing penetration of renewables
poses significant challenges to this model of operation because of its
intermittent and uncertain nature. At the same time, the penetration
of renewables also presents an opportunity to rethink the architecture
and its hierarchical separation towards the goal of improving
efficiency and adaptivity. A key aspect to achieve the integration of
different layers is the characterization of the robustness properties
of the mechanisms used at each layer, since variables at the upper
layers cannot be assumed in steady state any more at the lower ones.
These considerations motivate our work on iterative bidding schemes
combined with continuous physical network dynamics and the correctness
analysis of the resulting multi-rate hybrid interconnected system.



\subsubsection*{Literature review}

The integration of economic dispatch and frequency regulation in power
networks has attracted increasing attention in the last decades.  Many
recent
works~\cite{dorfler2016breaking,zhangpapaautomatica,trip2016internal,cady2015distributed,li2016connecting,zhang2017distributed,dall2016photovoltaic}
envision merging the design of primary, secondary, and tertiary
control layers for several models of the power network/micro-grid
dynamics with the aim of bridging the gap between long-term
optimization and real-time frequency control. In scenarios where
generators are price-takers, the literature has also explored the use
of market mechanisms to determine the optimal allocation of power
generation and to stabilize the frequency with real-time (locational
marginal) pricing,
see~\cite{alv_meng_power_coupl_market,DJS-MC-AMA:16,
  shiltz2017integration,stegink2017unifying}.  Inspired by the
iterative bidding schemes for strategic generators proposed
in~\cite{cherukuri2017iterative}, that lead to efficient Nash
equilibria where power generation levels minimize the total cost as
intended by the ISO, our work~\cite{TAC17-TS-AC-CDP-AS-JC} has shown
that the integration with the frequency dynamics of the network can
also be achieved in scenarios where generators are
price-bidders. However, this integration relies on a continuous-time
model for the bidding process, where the frequency coming from the
power network dynamics enters as a feedback signal in the negotiation
process. Instead, we account here for the necessarily discrete nature
of the bidding process and explore the design of provably correct
multi-rate hybrid implementations that realize this integration.


\subsubsection*{Statement of contributions}

We consider an electrical power network consisting of an ISO and a
group of strategic generators. The ISO seeks to ensure that the
generation meets the load with the minimum operation cost and the grid
frequency is regulated to its nominal value.  Each generator seeks to
maximize its individual profit and does not share its cost function
with anyone.  The ISO operates the market, where generators bid prices
at which there are willing to provide power, and makes power
generation assignments based on the bids and the local frequency
measurements.  Our goal is to design mechanisms that ensure the
stability of the interconnection between the ISO-generator bidding
process and the physical network dynamics while accounting for the
different nature (iterative in the first case, evolving in continuous
time in the second) of each process.  Our starting point is a
continuous-time bid update scheme coupled with the physical dynamics
of the power network whose equilibrium corresponds to an efficient
Nash equilibrium and zero frequency deviation.  Our first contribution
is the characterization of the robustness properties of this dynamics
against additive disturbances. To achieve this, we identify a novel
local Lyapunov function that includes the energy function of the
closed-loop system. The availability of this function not only leads
us to establish local exponential convergence to the desired
equilibrium, but also allows us rigorously establish its local
input-to-state stability properties.  Building on these results, our
second contribution develops a time-triggered hybrid implementation
that combines the discrete nature of iterative bidding with the
continuous nature of the frequency evolution in the power network.  In
our design, we introduce two iteration loops, one (faster) inner-loop
for the bidding process that incorporates at each step the frequency
measurements, and one (slower) outer-loop for the market clearing and
the updates in the power generation levels, that are sent to the
continuous-time power network dynamics.  We refer to this multi-rate
hybrid implementation as time-triggered because we do not necessarily
prescribe the time schedules to be periodic.  To analyze its
convergence properties, we regard the time-triggered implementation as
an approximation of the continuous-time dynamics and invoke the
robustness properties of the latter, interpreting as a disturbance
their mismatch. This allows us to derive explicit upper bounds on the
length between consecutive triggering times that guarantee that the
time-triggered implementation remains asymptotically convergent. The
computation of these upper bounds does not require knowledge of the
efficient Nash equilibrium.  Simulations on the IEEE 14-bus power
network illustrate our results.


\subsubsection*{Outline}
The paper is organized as follows.  Section~\ref{sec:phys-power-netw}
introduces the dynamic model of the power network and
Section~\ref{sec:problem-statement} describes the problem
setup. Section~\ref{sec:continuous-real-time} characterizes the
robustness properties of the continuous-time dynamics resulting from
the interconnection of bid updating and network dynamics.
Section~\ref{sec:iter-bidd-with} introduces the time-triggered
implementation and identifies sufficient conditions on the inter-event
times that ensure asymptotic convergence to efficient Nash
equilibria. Simulations illustrate the results in
Section~\ref{sec:simulations}. Section~\ref{sec:conclusions} gathers
our conclusions and ideas for future work. The appendices contain the
proofs of the main results of the paper.

\subsubsection*{Notation}
Let $\mathbb R , \mathbb R_{\ge0} , \mathbb R_{>0} , \mathbb Z_{\geq
  0} , \mathbb Z_{\ge 1}$ be the set of real, nonnegative real,
positive real, nonnegative integer, and positive integer numbers,
respectively.
For $m \in \mathbb Z_{\geq1}$, we use the shorthand notation
$\until{m} = \{1, \dots, m\}$.
For $A\in\mathbb R^{m\times n}$, we let $\norm{A}$ denote the induced
$2$-norm. Given $v\in\mathbb R^n,A=A^T\in\mathbb R^{n\times n}$, we
denote $\|v\|^2_A:=v^TA v$.  The notation $\1\in\mathbb R^n$ is used
for the vector whose elements are equal to 1.  The Hessian of a
twice-differentiable function $f:\mathbb R^n\to\mathbb R$ is denoted
by~$\n^2f$.
\section{Power network frequency
  dynamics}\label{sec:phys-power-netw}

Here we present the model of the physical power network that describes
the evolution of the grid frequency.  The network is represented by a
connected, undirected graph $ \G = (\V, \E) $, where nodes $ \V =
\until{n}$ represent buses and edges $ \E \subset \V \times \V $ are
the transmission lines connecting the buses.  Let $m$ denote the
number of edges, arbitrarily labeled with a unique identifier in
$\until{m}$.
The ends of each edge are also arbitrary labeled with ‘+’ and ‘-’, so
that we can associate to the graph the incidence matrix $D \in
\real^{n \times m}$ given by
\begin{align*}
  D_{ik}=
  \begin{cases}
    +1 &\text{if $i$ is the positive end of edge $k$},\\
    -1 &\text{if $i$ is the negative end of edge $k$},\\
    0 & \text{otherwise.}
  \end{cases}
\end{align*}
Each bus represents a control area and is assumed to have one
generator and one load. Following~\cite{powsysdynwiley}, the dynamics
at the buses is described by the \emph{swing equations}  \eqref{eq:swingdeltacomp}.
\begin{equation}
  \begin{aligned}
    \dot \delta&=\w\\
    M\dot \w&= -D\Gamma\sin (D^T\delta)-A\w+P_g-P_d
  \end{aligned}\label{eq:swingdeltacomp}
\end{equation}
Here $\Gamma=\diag\{\gamma_1,\ldots,\gamma_m\}\in\mathbb R^{m\times
  m},\gamma_k=B_{ij}V_{i}V_{j}$, where $k\in\until{m}$ corresponds to
the edge between nodes $i$ and~$j$.  Table~\ref{tab:par3SG} specifies
the meaning of the symbols used in the
model~\eqref{eq:swingdeltacomp}.
\begin{table}[htb]\vspace*{-3mm}
  \begin{align*}
    \delta&\in \mathbb R^n && \text{(vector of) voltage phase angles}  \\
    \w&\in\mathbb R^n & &  \text{frequency deviation w.r.t. the nominal frequency}  \\
    V_{i}&\in\mathbb R_{>0}& & \text{voltage magnitude at bus } i
    \\ 
    P_{d}&\in\mathbb R^n  && \text{power load}  \\
    P_{g}&\in\mathbb R^n & & \text{power generation} \\
    M&\in\mathbb R_{\ge0}^{n\times n} & & \text{diagonal matrix of moments of inertia} \\
    A&\in\mathbb R_{\geq0}^{n\times n}&  & \text{diagonal matrix of asynchronous damping constants} \\
    B_{ij}&\in\mathbb R_{\geq0}& & \text{negative of the susceptance of
      transmission line } (i,j)
  \end{align*}
  \caption{Parameters and state variables of model~\eqref{eq:swingdeltacomp}.}\label{tab:par3SG}
\end{table}

To avoid issues in the stability analysis of \eqref{eq:swingdeltacomp}
due to the rotational invariance of $\delta$, see
e.g.,~\cite{dorfler2014synchronization}, we introduce the new variable
$\vp=D_t^T\delta\in\mathbb R^{n-1}$. Here $\vp$ represents the voltage
phase angle differences along the edges of a spanning tree of the
graph $\G$ with incidence matrix $D_t$.  The physical energy stored in
the transmission lines is given by~\eqref{eq:U}, where
$D_t^\dagger=(D_t^TD_t)^{-1}D_t^T$ denotes the Moore-Penrose inverse
of $D_t$.
\begin{align}\label{eq:U}
  U(\vp)=-\1^T\Gamma\cos(D^TD_t^{\dagger T}\vp) .
\end{align}
By noting that $D_tD_t^{\dagger}D=(I-\frac1n\1\1^T)D=D$, the physical
system~\eqref{eq:swingdeltacomp} in the $(\vp,\w)$-coordinates takes
the form
\begin{equation}\label{eq:swingeqU}
  \begin{aligned}
    \dot \vp&=D_t^T\w\\
    M\dot \w&=-D_t\nabla U(\vp)-A\w+P_g-P_d
  \end{aligned}
\end{equation}
In the sequel we assume that, for the power generation $P_g=\bar P_g$,
there exists an equilibrium $\col(\bar \vp,\bar \w)$
of~\eqref{eq:swingeqU} that satisfies $D^TD_t^{\dagger
  T}\bar\vp\in(-\frac\pi2,\frac\pi2)^m$. The latter assumption is
standard and often referred to as the \emph{security
  constraint} 
\cite{powsysdynwiley}. 

\section{Problem statement}\label{sec:problem-statement}

In this section we formulate the problem statement and then discuss
the paper objectives. We start from the power network model introduced
in Section~\ref{sec:phys-power-netw} and then explain the
game-theoretic model describing the interaction between the ISO and
the generators following the exposition
of~\cite{cherukuri2016decentralized,cherukuri2017iterative}.

The cost incurred by generator $i \in [n]$ in producing $P_{gi}$ units
of power is given by
\begin{align}\label{eq:cost-Ci}
  C_i(P_{gi}):=\frac12 q_i P_{gi}^2 + c_i P_{gi} ,
\end{align}
where $q_i >0$ and $c_i \ge 0$. The total network cost is then
\begin{align}\label{eq:total-cost}
  C(P_g) :=\sum_{i\in[n]}C_i(P_{gi})=\frac12P_g^TQP_g+c^TP_g,
\end{align}
with $Q=\diag\{q_1,\ldots,q_n\}$ and $c=\col(c_1, \dots, c_n)$. Given
the cost~\eqref{eq:total-cost} and the power loads $P_d$, the ISO
seeks to solve the \emph{economic dispatch problem}
\begin{subequations}\label{eq:ISO-OPFsc}
  \begin{align}
    \underset{P_g}{\minimize} & \quad C(P_g),
    \\
    \st & \quad \1^TP_g=\1^TP_d ,
    \label{eq:ISO-nodal-pow-bal}
  \end{align}
\end{subequations}
and, at the same time, regulate the network frequency to its nominal
value.  Since the function $C$ is strongly convex, there exists a
unique optimizer $P_g^*$ of~\eqref{eq:ISO-OPFsc}. However, we assume
that the generators are strategic and they do not reveal their cost
functions to anyone, including the ISO. Consequently, the ISO is
unable to determine the optimizer of~\eqref{eq:ISO-OPFsc}.  Instead,
it determines the power dispatch according to a market clearing
procedure in which each generator submits bids to the ISO. 

We consider price-based bidding: each generator $i \in \until{n}$
submits the price per unit electricity $b_i \in \real$ at which it is
willing to provide power.  Based on these bids, the ISO finds the
power generation allocation that minimizes the total generator payment
while meeting the load. More precisely, given the bid
$b=\col(b_1,\ldots,b_n)$, the ISO solves
%
%
%
%
\begin{subequations}\label{eq:ISOprob}
  \begin{align}
    \underset{P_g}{\minimize} & \quad b^TP_g,  \label{eq:ISOobj-fun} 
    \\
    \st & \quad \1^TP_g=\1^TP_d. \label{eq:pow-bal-con}
  \end{align}
\end{subequations}
The optimization problem \eqref{eq:ISOprob} is linear and may in
general have multiple (unbounded) solutions.
Among these solutions, let
$P_{g}^{\text{opt}}(b)=\col(P_{g1}^{\text{opt}}(b),\ldots,P_{gn}^{\text{opt}}(b))$
be the optimizer of \eqref{eq:ISOprob} the ISO selects given
bids~$b$. Knowing this process, each generator $i$ aims to bid a
quantity $b_i$ to maximize its payoff
\begin{equation}\label{eq:payoffgen}
  \Pi_i(b_i,P_{gi}^{\text{opt}}(b)) :=P_{gi}^{\text{opt}}(b)b_i-C_i(P_{gi}^{\text{opt}}(b)).
\end{equation}
For an unbounded optimizer we have $\Pi_i(b_i,\pm \infty)=-\infty$.
To analyze the clearing of the market, we resort to tools from game
theory~\cite{
  DF-JT:91}. To this end, we define the \emph{inelastic electricity
  market game}:
\begin{itemize}
\item Players: the set of generators $[n]$.
\item Action: for each player $i\in[n]$, the bid $b_i\in\mathbb R$.
\item Payoff: for each player $i\in[n]$, the payoff $\Pi_i$ in
  \eqref{eq:payoffgen}.
\end{itemize}
For the bid vector we interchangeably use the notation $b\in\mathbb
R^n$ and $(b_i,b_{-i})\in\mathbb R^n$, where $b_{-i}$ represents the
bids of all players except $i$. 
%
A bid profile $b^*\in\mathbb R^n$ is a \emph{Nash equilibrium} if
there exists an optimizer $P_{g}^{\text{opt}}(b^*)$
of~\eqref{eq:ISOprob} such that $\forall i\in[n]$,
\begin{align*}
  \Pi_i(b_i,P_{gi}^{\text{opt}}(b_i,b_{-i}^*))\leq
  \Pi_i(b_i^*,P_{gi}^{\text{opt}}(b^*))
\end{align*}
for all $b_i\neq b_i^*$ and all optimizers
$P_{gi}^{\text{opt}}(b_i,b_{-i}^*)$ of \eqref{eq:ISOprob}. 
%
In particular, we are interested in bid profiles that can be
associated to economic dispatch.
More specifically, a bid $b^*\in\mathbb R^n$ is \emph{efficient} is a
bid if there exists an optimizer $P_g^*$ of \eqref{eq:ISO-OPFsc} which
is also an optimizer of \eqref{eq:ISOprob} given bids $b=b^*$ and
\begin{align}\label{eq:eff-bid}
  P_{gi}^*=\argmax_{P_{gi}}\{P_{gi} b^*_i-C_i(P_{gi})\} \text{ for all } i
  \in [n].
\end{align} 
A bid $b^*$ is an \emph{efficient Nash equilibrium} if it is both
efficient and a Nash equilibrium. At the efficient Nash equilibrium,
the optimal generation allocation determined by \eqref{eq:ISO-OPFsc}
coincides with the production that the generators are willing to
provide, maximizing their profit~\eqref{eq:payoffgen}.  Following the
same arguments as in the proof of \cite[Lemma
3.2]{cherukuri2016decentralized}, 
one can establish the existence and uniqueness of the efficient Nash
equilibrium.
\begin{proposition}\longthmtitle{Existence and uniqueness of efficient
    Nash equilibrium}\label{prop:exis-NE}
  Let $( P_{g}^*,\l^*)$ be a primal-dual optimizer of
  \eqref{eq:ISO-OPFsc}, then $b^* \! = \! \nabla C(P_g^*) \! = \!
  \1\l^*$ is the unique efficient Nash equilibrium of the inelastic
  electricity market game.
\end{proposition}



In the scenario described above, neither the ISO nor the individual
strategic generators are able to determine the efficient Nash
equilibrium beforehand.  Our goal is then to design an online bidding
algorithm where ISO and generators iteratively exchange information
about the bids and the generation quantities before the market is
cleared and dispatch commands are sent. The algorithm should be truly
implementable, meaning that it should account for the discrete nature
of the bidding process, and at the same time ensure that network
frequency, governed by the continuous-time power system dynamics, is
regulated to its nominal value. The combination of these two aspects
leads us to adopt a hybrid implementation strategy to tackle the
problem.

\section{Robustness of the continuous-time bid and power-setpoint
  update scheme}\label{sec:continuous-real-time}

In this section, we introduce a continuous-time dynamics that
prescribes a policy for bid updates paired with the frequency dynamics
of the power network whose equilibrium corresponds to an efficient
Nash equilibrium and zero frequency deviation. In this scheme,
generators update their bids in a decentralized fashion based on the
power generation quantities received by the ISO, while the ISO changes
the generation quantities depending on both the generator bids and the
network frequency.  This design is a simplified version of the one
proposed in our previous work~\cite{TAC17-TS-AC-CDP-AS-JC}. The main
contribution of our treatment here is the identification of a novel
Lyapunov function that, beyond helping establish local exponential
convergence, allows us to characterize the  input-to-state
stability properties of the dynamics. We build on this
characterization later to develop our time-triggered hybrid
implementation that solves the problem outlined in
Section~\ref{sec:problem-statement}.

\subsection{Bidding process coupled with physical network dynamics}

Recall from Section~\ref{sec:problem-statement} that given bid~$b_i$,
generator $i\in[n]$ wants to produce the amount of power that
maximizes its individual profit, given by
\begin{align}
  P_{gi}^{\text{des}}:=\argmax_{P_{gi}}\{b_iP_{gi}-C_i(P_{gi})\}
  =q_i^{-1}(b_i-c_i)\label{eq:Pg-des}
\end{align}
Hence, if the ISO wants generator $i$ to produce more power than its
desired quantity, that is $P_{gi}>P_{gi}^{\text{des}}$, generator $i$
will increase its bid, and vice versa. Bearing this rationale in mind,
the generators update their bids according to
\begin{subequations}\label{eq:bidding-process}
  \begin{align}\label{eq:produpdateb}
    \gain_b \dot b&=P_g-Q^{-1}b+Q^{-1}c.
  \end{align}
  Here $\gain_b\in\mathbb R^{n\times n}$ is a diagonal positive
  definite matrix.
  Next, we provide an update law for the ISO depending on the bid
  $b\in\mathbb R^n$ and the local frequency of the power network. 
  The ISO updates its actions according to
  \begin{equation}
    \begin{aligned}
      \gain_g\dot P_g&=\1\l-b+\rho \1\1^T(P_d-P_g) -\sigma^2 \w\\
      \gainl_\l\dot \l&=\1^T(P_d-P_g)
    \end{aligned}\label{eq:ISOprimaldual}
  \end{equation}
\end{subequations}
with parameters $\rho,\sigma,\gainl_\l\in\mathbb R_{>0}$ and where
$\gain_g\in\mathbb R^{n\times n}$ is a diagonal positive definite gain
matrix.

The intuition behind the dynamics \eqref{eq:ISOprimaldual} is
explained as follows. If generator $i$ bids higher than the 
Lagrange multiplier $\l$ (sometimes referred to as the \emph{shadow
  price} \cite{stoft2002power}) associated to \eqref{eq:pow-bal-con},
then the power generation (setpoint) of node $i$ is decreased, and
vice versa. By adding the term with $\rho>0$, one can enhance the
convergence rate of \eqref{eq:ISOprimaldual}, see
e.g.,~\cite{boyd2011distributed}.
We 
add the feedback signal $-\sigma^2\w$ to compensate for the frequency
deviations in the physical
system. 
Interestingly, albeit we do not pursue this here, the dynamics
\eqref{eq:bidding-process} could also be implemented in a distributed
way without the involvement of a central regulating authority like the
ISO.

For the remainder of the paper, we assume that there exists an
equilibrium $\bar x=\col(\bar \vp,\bar \w,\bar b,\bar P_g,\bar \l)$
of~\eqref{eq:swingeqU}-\eqref{eq:bidding-process} such that
$D^TD_t^{\dagger T}\bar\vp\in(-\pi/2,\pi/2)^m$ (cf.
Section~\ref{sec:phys-power-netw}). Note that this equilibrium
satisfies
\begin{align}\label{eq:eq-cont-time}
  \begin{aligned}
    \bar\l&=\frac{\1^T(P_d+Q^{-1}c)}{\1^TQ^{-1}\1}>0, &\bar \w&=0,
    \quad \bar b=\1\bar \l,
    \\
    \bar P_g&=Q^{-1}\1\bar\l-Q^{-1}c, & \1^T\bar P_g&=\1^T P_d.
  \end{aligned}
\end{align}
In particular, at the steady state, the frequency deviation is zero,
the power balance $\1^T\bar P_g=\1^TP_d$ is satisfied, and
$\1\bar \l=\bar b=\nabla C(\bar P_g)$, implying that $\bar P_g$ is a
primal optimizer of \eqref{eq:ISO-OPFsc} and $\bar b$ is an efficient
Nash equilibrium by Proposition \ref{prop:exis-NE}. Hence, at steady
state the generators do not have any incentive to deviate from the
equilibrium bid.

\subsection{Local input-to-state (LISS) stability}

While the ISO dynamics \eqref{eq:ISOprimaldual} is a saddle-point
dynamics of the linear optimization problem \eqref{eq:ISOprob} (and
hence, potentially unstable), we show next that the interconnection of
the physical power network dynamics~\eqref{eq:swingeqU} with the
bidding process~\eqref{eq:bidding-process} is locally exponentially
stable and, furthermore, robust to additive disturbances.  For
$x=\col(\vp,\w,b,P_g,\l)$, define the function
\begin{align}\label{eq:defV}
  V(x)
  &=U(\vp)-(\vp-\bar \vp)^T\n U(\bar \vp)-U(\bar \vp)+\tfrac12\w^TM\w\nonumber
  \\
  &+\tfrac1{2\sigma^2}(\|b-\bar b\|^2_{\gain_b}+
    \|P_g-\bar P_g\|^2_{\gain_g}+
    \|\l-\bar \l\|_{\gainl_\l}^2).
\end{align}
Then the closed-loop system obtained by combining~\eqref{eq:swingeqU}
and~\eqref{eq:bidding-process} is compactly written as
\begin{align}\label{eq:compactclsysdyn}
  \dot{x}&=F(x)= \mathcal Q^{-1}
  \mathcal A \mathcal Q^{-1}\nabla V(x)
\end{align}
%
with
$\mathcal Q=\mathcal
Q^T=\blockdiag(I,M,\tfrac{\gain_b}\sigma,\tfrac{\gain_g}\sigma,\tfrac{\gainl_\l}\sigma)>0$
and
\begin{align*}
  \mathcal A=
  \begin{bmatrix}
    0     & D_t^T      & 0    & 0        & 0 \\
    - D_t & - A        & 0    & \sigma I & 0 \\
    0     & 0          & - Q^{-1}  & I        & 0 \\
    0     & - \sigma I & -  I & -\rho \1\1^T  & \1 \\
    0     & 0          & 0    & - \1^T      & 0
  \end{bmatrix}.
\end{align*}
By exploiting the structure of the system, we obtain the dissipation
inequality
\begin{align}\label{eq:V-diss-ineq}
  \dot V=\frac12(\nabla V(x))^T\mathcal Q^{-1}(\mathcal A+\mathcal
  A^T)\mathcal Q^{-1}\nabla V(x)\leq 0
\end{align}
However, since $\mathcal R:=-\frac12(\mathcal A+\mathcal A^T)$ is only
positive semi-definite, $V$ is not strictly decreasing along the
trajectories of~\eqref{eq:compactclsysdyn}.  Nevertheless, one can
invoke the LaSalle Invariance Principle to characterize the local
asymptotic convergence properties of the dynamics,
cf.~\cite{TAC17-TS-AC-CDP-AS-JC}. Here, we show that, in fact, the
dynamics is locally input-to-state (LISS) stable, as defined
in~\cite{EDS-YW:96}, and therefore robust to additive
disturbances. Our key tool to establish this is the identification of
a LISS-Lyapunov function, which in general is far from trivial for
dynamics that involve saddle-point dynamics. To this end, consider the
system
\begin{align}\label{eq:compactclsysdyn-dist}
  \dot{x}&=F(x)+Bd
\end{align}
with $B\in\real^{4n\times q}$ and a disturbance signal
$d\in\real^q$. In the following result, we use the function $V$ to
construct an LISS-Lyapunov function for the system~\eqref{eq:compactclsysdyn-dist}.


\begin{theorem}\longthmtitle{LISS-Lyapunov function for the
    interconnected dynamics}\label{thm:VeLyap}
  Consider the interconnected dynamics~\eqref{eq:compactclsysdyn-dist}
  and define the function
  \begin{align}
    &W_\e(x) = V(x)+\e_0\e_1(\vp-\bar
    \vp)^TD_t^{\dagger}M\w\label{eq:def-We}\\
    &-\tfrac{\e_0\e_2}{\sigma^2}
    (b-\bar b)^T \gain_g(P_g-\bar P_g)
    - \tfrac{\e_0\e_3}{\sigma^2}
    (\l-\bar \l)\1^T \gain_g(P_g-\bar P_g) , \nonumber
  \end{align}
  with parameters $\e=\col(\e_0,\e_1,\e_2,\e_3)\in\mathbb R^4_{>0}$
  and $V$ given by \eqref{eq:defV}.  Given the equilibrium
  $\bar x=\col(\bar \vp,\bar \w,\bar b,\bar P_g,\bar \l)$ of
  \eqref{eq:compactclsysdyn}, let
  $\bar \eta=D^TD_t^{\dagger T}\bar\vp$. For $\gamma$ such that
  $\norminf{\bar \eta}<\gamma<\frac\pi2$, define the closed convex set
  \begin{align}\label{eq:def-Omega}
    \Omega=\{x=\col(\vp,\w,b,P_g,\l) \ \! |\! \ D^TD_t^{\dagger
      T}\!\vp \in [-\gamma,\gamma ]^m\}.
  \end{align}
  Then there exist sufficiently small $\e$ such that $W_\e$ is an
  LISS-Lyapunov function of \eqref{eq:compactclsysdyn-dist}
  on~$\Omega$. In particular, there exist constants $\alpha,\chi,
  c_1,c_2>0$ such that for all $x\in \Omega$ and all $d$ satisfying
  $\norm d\leq \chi\norm {x-\bar x}$,
  \begin{subequations}\label{eq:q-Lyapunov}
    \begin{align}
      \tfrac12 c_1\|x-\bar x\|^2\leq W_\e(x)&\le \tfrac12 c_2\|x-\bar
      x\|^2, \label{eq:quadWe}
      \\
      (\nabla W_\e(x))^T (F(x)+Bd)&\leq -\alpha \|x-\bar
      x\|^2
      . \label{eq:xiein}
    \end{align}
  \end{subequations}
\end{theorem}

We refer to Appendix \ref{sec:proof-theor-refthm:v-1} for the proof of
Theorem~\ref{thm:VeLyap}.  Using the
characterization~\eqref{eq:q-Lyapunov} and \cite[Theorem
4.10]{khalil1996nonlinearLaSalle}, each trajectory of
\eqref{eq:compactclsysdyn} initialized in a compact level set
contained in $\Omega$ exponentially converges to the equilibrium $\bar
x$ corresponding to economic dispatch and the efficient Nash
equilibrium. Moreover, we exploit the local ISS property
of~\eqref{eq:compactclsysdyn-dist} guaranteed by
Theorem~\ref{thm:VeLyap} next to develop a time-triggered hybrid
implementation.

\section{Time-triggered implementation: iterative bid update and
  market clearing}\label{sec:iter-bidd-with}

In realistic implementations, the bidding process between the ISO and
the generators is not performed continuously. Given the availability
of digital communications, it is reasonable to instead model it as an
iterative process.  Building on the continuous-time bidding dynamics
proposed in Section~\ref{sec:continuous-real-time}, here we develop a
time-triggered hybrid implementation that combines the discrete nature
of bidding with the continuous nature of the frequency evolution in
the power network. We consider two time-scales, one (faster) for the
bidding process that incorporates at each step the frequency
measurements, and another one (slower) for the market clearing and
updates of the power generation levels that are sent to the power
network dynamics.
We refer to this implementation as time-triggered because we do not
necessarily prescribe the time schedules to be periodic in order to
guarantee that the asymptotic stability properties are retained by the
hybrid implementation.

\subsection{Algorithm description}

We start with an informal description of the iterative update scheme
between the ISO and the generators, and the interconnection with the
dynamics of the power network.


\begin{quote}
  \emph{[Informal description]:} The algorithm has two time indices,
  $k$ to label the iterations on the bidding process and $l$ to label
  the iteration in the market clearing process that updates the power
  setpoints.  At each iteration $l\in \mathbb Z_{\ge0}$, ISO and
  generators are involved in an iterative process where, at each
  subiteration~$k$, generators send a bid to the ISO. Once the ISO has
  obtained the bids and the network frequency measurements at time
  $t_{k}^l$, it computes the new potential generation allocations,
  denoted $P_g^{k+1}\in\mathbb R^n$, and sends the corresponding one
  to each generator. At the $(k+1)$-th subiteration, each generator
  adjust its bid based on their previous bid and the generation
  allocation received from the ISO at time $t_{k+1}^l$. Once
  $k=N_l\in\mathbb Z_{\geq 1}$ at time $t_{N_l}^l$, the market is
  cleared, meaning that the bidding process is reset (i.e., $k=0$),
  the power generations in the swing equations are updated according
  to the current setpoints $P_g^{N_l}$, and the index~$l$ moves
  to~$l+1$.
\end{quote}
%
\begin{figure*}[htb]
  \centering
  \begin{tikzpicture}
    \foreach \x/\k/\l in
    {1/1/0,1.8/2/0,3.2/k-1/0,3.85/k/0,4.6/k+1/0,7.2/1/1,8.3/k/l} {
      \draw (\x ,0) node[anchor=south]{$t_{\k}^\l$} -- (\x,-0.1);
      \draw[<-] (\x ,-0.3) -- (\x,-1.9); };
    \draw[very thick] (6 ,0) node[anchor=south]
    {$t_{N_{0}}^{0}=t_0^{1}$} -- (6,-0.2);
    \draw[very thick,<->] (6 ,-0.3) -- (6,-1.8);
    \draw[very thick] (10 ,0) node[anchor=south]
    {$t_{N_{l}}^{l}=t_0^{l+1}$} -- (10,-0.2);
    \draw[very thick,<->] (10 ,-0.3) -- (10,-1.8);
    \draw[very thick] (0 ,0) node[anchor=south] {$t_{0}^0=0$} --
    (0,-0.2);
    \draw[very thick] (0 ,-2) -- (0,-1.9);
    \draw[very thick] (6 ,-2) -- (6,-1.9);
    \draw[very thick] (10 ,-2) -- (10,-1.9);
    \draw[very thick,<->] (0 ,-0.3)  -- (0,-1.8);
    \draw[-] (-.5,0) node(0){};
    \foreach \n/\i/\j in {2.4/0/1,5.2/1/2,7.5/2/3,8.7/3/4}{
      \draw[-] (\i) -- (\n,0) node[right](\j){$.$};
    };
    \draw[->] (4) -- (11,0);
    \draw[->] (-.5,-2) node[below right]{power network dynamics
      \eqref{eq:swingeqU}} -- (11,-2) node[below left]{$t$};
    \draw (-0.5,0.6) node[above right]{ISO-generator bidding process
      (Algorithm \ref{alg:bid-adjustm-algor})};
    \node(ISO) at (-3.5,0.6) [thick,draw,minimum width=1cm,minimum
    height=1cm,text width=2.5cm]{ISO-generator bidding process};
    \node(freq) at (-3.5,-1.9) [thick,draw,minimum width=1cm,minimum
    height=1cm,text width=2.5cm]{power network dynamics};
    \draw[->,transform canvas={xshift=-5mm}] (freq) --node[left,text
    width =1.5cm]{frequency deviations} (ISO);
    \draw[very thick,->,transform canvas={xshift=5mm}] (ISO)
    --node[right,text width =1.5cm]{power generation setpoints}
    (freq);
  \end{tikzpicture}
  \caption{Relation between time and iteration numbers in the
    time-triggered system \eqref{eq:timetrigclsys}. The lower
    time-axis corresponds to the continuous-time physical system
    \eqref{eq:swingeqU} while the upper one corresponds to the time
    sequence $\{\{t_k^l\}_{k=0}^{N_l}\}_{l=0}^\infty$ of the
    ISO-generator bidding process given in Algorithm
    \ref{alg:bid-adjustm-algor}. The arrows pointing up indicate the
    frequency updates in the bidding dynamics while the arrows
    pointing down correspond to update of the power generation levels
    in the physical system. As indicated, for each
    $l\in\mathbb Z_{\ge0}$ the lower index $k$ is reset once it
    reaches $k=N_l\in\mathbb Z_{\ge1}$, i.e.,
    $t_{N_{l}}^{l}=t_0^{l+1}$ for all $l\in\mathbb
    Z_{\ge0}$. }\label{fig:time-instances-hybrid}
\end{figure*}
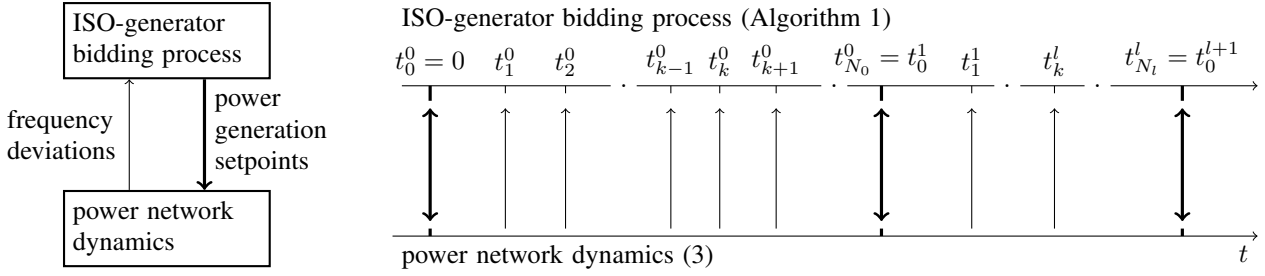

Figure~\ref{fig:time-instances-hybrid} shows the two iteration layers
in the update scheme. The evolution of the frequency occurs in
continuous time according to~\eqref{eq:swingeqU}. To relate iteration
numbers with time instances on $\mathbb R$, we consider time sequences
of the form $\{ \{t_k^l\}_{k=0}^{N_l} \}_{l=0}^\infty$ for
$N_l\in\mathbb Z_{\ge1}$ and $l\in\mathbb Z_{\ge0}\}$, satisfying
\begin{align}\label{eq:time-seq-ord}
  t_{k}^l-t_{k-1}^l>0, \qquad t_{0}^{l+1}=t_{N_{l}}^{l} 
\end{align}
for all $ l\in\mathbb Z_{\ge0}$ and all $ k\in[N_l]$.
Algorithm~\ref{alg:bid-adjustm-algor} formally describes the iterative
updates of the bidding process between the generators and the ISO.

\removelatexerror 
\begin{algorithm}[h]
  \SetAlgoLined
  \DontPrintSemicolon
  \SetKwFor{Case}{case}{}{endcase}
  \SetKwInOut{giv}{Data} \SetKwInOut{ini}{Initialize}
  \SetKwInOut{state}{State} \SetKwInput{start}{Initiate}
  \SetKwInOut{msg}{Messages} \SetKwInput{Kw}{Executed by}
  \Kw{generators $i \in \until{n}$ and ISO}
  \giv{time sequence $\{ \{t_k^l\}_{k=0}^{N_l} \}_{l=0}^\infty $;
    cost function \eqref{eq:cost-Ci}
    for each generator $i$; frequency deviation
    $\w(t_k^l)$ at each time $t_k^l$ and load $P_d$ for ISO}
  %
  \ini{each generator $i$ selects arbitrarily $b_i^0 \ge \coeffc_i$,
    sets $k=0,l=0$, and jumps to step~\ref{step1}; ISO selects
    arbitrary $P_{gi}^0>0,\l_i^0>0$, sets $k=0,l=0$ and waits for
    step~\ref{step2}}
  \BlankLine \While{$l\geq0$}{
    \BlankLine \While{$k\geq0, k< N_l$}{
      \tcc{For each generator $i$:} Receive $P_{gi}^k$ from ISO at
      $t_k^l$; Set\;
      $b_i^{k+1}=b_i^k+(t_{k+1}^l-t_k^l)\gain_{bi}^{-1}(P_{gi}^k-q_i^{-1}(b_i^k+c_i))$ \label{step3}\;
      %
      %
      Send $b_i^{k+1}$ to the ISO; set $k = k+1$ \label{step1}\;
      \BlankLine \tcc{For ISO:} Receive $b_i^k,\w_i(t_k^l)$ from each
      $i \in \until{n}$ at $t_k^l$ \;
      \label{step2} Set
      $P_{gi}^{k+1}=P_{gi}^k+(t_{k+1}^l-t_k^l)\gain_{gi}^{-1}(\l^k-b_i^k-\sigma^2\w_i(t_k^l)
      +\rho \mysum[]_{i\in[n]}(P_{di}-P_{gi}^k))$ for all $ i\in[n]$
      $\l^{k+1}= \l^k+\frac{t_{k+1}^l-t_k^l}{\gainl_\l}
      \mysum[]_{i\in[n]}(P_{di}-P_{gi}^k)$ \;
      Send $P_{gi}^{k+1}$ to each $i \in \until{n}$, set $k=k+1$ \; }
    Set $P_{gi}(t)=P_{gi}^{N_l}$ in \eqref{eq:swingeqU}
    $\forall i\in[n], \forall t\in [t_{N_l}^l,t_{N_{l+1}}^{l+1})$\;
    Set $b_i^0=b_i^{N_l},P_{gi}^0=P_{gi}^{N_l}, \l_{i}^0=\l_{i}^{N_l}$
    for each $i\in [n]$\; Set $l=l+1, k=0$\;} 
  \caption{\baalgo}\label{alg:bid-adjustm-algor}
\end{algorithm}
\vspace*{1ex}

%
%
%
%
%
%
%
For analysis purposes, we find it convenient to represent the dynamics
resulting from the combination of
Algorithm~\ref{alg:bid-adjustm-algor} and the network
dynamics~\eqref{eq:swingeqU} as the time-triggered continuous-time
system
\begin{align}
  \dot \vp(t)
  &=D_t^T\w(t),\nonumber
  \\
  M\dot \w(t)
  &=-D_t\nabla U(\vp(t))-A\w(t)+P_g(t_{0}^l)- P_d,\nonumber
  \\
  \gain_b\dot
  b(t)
  &=P_g(t_k^l)-Q^{-1}b(t_k^l)-Q^{-1}c,\label{eq:timetrigclsys}
  \\
  \gain_g\dot P_g(t)
  &=\1\l(t_k^l)-b(t_k^l)-\sigma^2 \w(t_k^l)  
        +\rho \1\1^T(P_d-P_g(t_k^l)),\nonumber
  \\
  \gainl_\l\dot \l(t)
  &=
    \1^T(P_d-P_g(t_k^l)),
    \nonumber
\end{align}
for
$t\in [t_{k}^l,t_{k+1}^l)\subset [t_0^l,t_{0}^{l+1}), l\in\mathbb
Z_{\ge0}, k\in \{0,\ldots,N_l-1\}$.
%
%
We write the system~\eqref{eq:timetrigclsys} compactly in the form
\begin{align}\label{eq:gentwotrigsys}
  \dot x(t)&=f(x(t))+g(x(t_k^l))+h(x(t_0^l))
\end{align}
with
\begin{align*}
  & f(x)=\col(D_t^T\w,-M^{-1}(D_t\nabla U(\vp)+A\w+P_d),0,0,0),
  \\
  &g(x)=\col(0,0,\gain_b^{-1}(P_g-Q^{-1}b-Q^{-1}c),
  \\
  &\gain_g^{-1}( \1 \l-b-\sigma^2\w +
    \rho\1\1^T(P_d-P_g)),\gainl_\l^{-1}\1^T(P_d-P_g)),
  \\
  & h(x)=\col(0,M^{-1}P_g,0,0,0).
\end{align*}

With this notation, note that the continuous-time
dynamics~\eqref{eq:compactclsysdyn} corresponds to
\begin{align}\label{eq:gennonsys}
  \dot x(t)= f(x(t))+g(x(t))+h(x(t)).
\end{align}
Since $\sup_{\vp\in\mathbb R^{n-1}}\|\nabla^2U(\vp)\|<\infty$
and $g,h$ are linear, it follows that $f,g,h$ are globally
Lipschitz (we denote by $L_f,L_g,L_h$ their Lipschitz constants,
respectively).
When viewed as a continuous-time system, the
dynamics~\eqref{eq:timetrigclsys} has a discontinuous right-hand side,
and therefore we consider its solutions in the Carath\'eodory sense,
cf.~\cite{JC:08-csm}.



\subsection{Sufficient condition on triggering times for stability}

In this section we establish conditions on the time sequence that
guarantee that the solutions of~\eqref{eq:timetrigclsys} are
well-defined and retain the convergence properties
of~\eqref{eq:compactclsysdyn}.  Specifically, we determine a
sufficient condition on the inter-sampling times $t_{k+1}^l- t_k^l$
for bidding and $t_k^{l+1} - t_k^l$ for market clearing that ensure
local asymptotic convergence of~\eqref{eq:gentwotrigsys} to the
equilibrium~$\bar x$ of the continuous-time
system~\eqref{eq:compactclsysdyn}.

Our strategy to accomplish this relies on the robustness properties
of~\eqref{eq:compactclsysdyn} characterized in
Theorem~\ref{thm:VeLyap} and the fact that the time-triggered
implementation, represented by~\eqref{eq:gentwotrigsys}, can be
regarded as an approximation of the continuous-time dynamics,
represented by~\eqref{eq:gennonsys}. We use the Lyapunov function
$W_\e$ defined by \eqref{eq:def-We} and examine the mismatch between
both dynamics to derive upper bounds on the inter-event times that
guarantee that $W_\e$ is strictly decreasing along the time-triggered
system~\eqref{eq:timetrigclsys}.

\begin{theorem}\longthmtitle{Local asymptotic stability of
    time-triggered implementation}\label{thm:main-conv-result}
  Consider the time-triggered implementation~\eqref{eq:timetrigclsys}
  of the interconnection between the ISO-generator bidding processes
  and the power network dynamics. With the notation of
  Theorem~\ref{thm:VeLyap}, let
  \begin{align}
    \bar \xi &:= \frac{1}{L_f + L_g} \log \Bigl(1+ \frac{\beta(L_f +
               L_g)}{L(L_W L_h + \beta)} \Bigr), \label{eq:xibar-zetabar-def} 
    \\
    \bar \zeta &:= \frac{1}{L_f} \log \Bigl( 1+ 
                 \frac{ L_f (\alpha
                 - \beta)}{L_g (L L_W + \alpha) + (\alpha - \beta) (L_f +
                 L_g)} \Bigr), \notag 
  \end{align}
  where $0<\beta<\alpha$, $L:=L_f+L_g+L_h$, and $L_W$ is the Lipschitz
  constant of $\nabla W_\e$.  Assume the time sequence
  $\{ \{t_k^l\}_{k=0}^{N_l} \}_{l=0}^\infty$ satisfies, for some
  $\zetaubar \in (0,\zetaobar)$ and $\xiubar \in (0,\xiobar)$,
  \begin{align}\label{eq:step-size-bounds}
    \zetaubar \le t_0^{l+1} - t_0^l \le \zetaobar \quad \text{and} \quad
    \xiubar \le t_{k}^l - t_{k-1}^l \le \xiobar,
  \end{align}
  for all $l \in \mathbb Z_{\ge0}$ and $k \in \until{N_l}$.
  Then, $\bar x$ is locally asymptotically stable
  under~\eqref{eq:timetrigclsys}.
\end{theorem}

We refer the reader to Appendix~\ref{sec:two-time-triggering} for the
proof of Theorem~\ref{thm:main-conv-result}.  The uniform lower bounds
$\ul{\zeta}$ and $\ul{\xi}$ in~\eqref{eq:step-size-bounds} ensure that
the solutions of the time-triggered implementation
\eqref{eq:timetrigclsys} are well-defined, avoiding Zeno behavior.
Theorem~\ref{thm:main-conv-result} implies that convergence is
guaranteed for any constant stepsize implementation, where the
sufficiently small stepsize
satisfies~\eqref{eq:step-size-bounds}. However, the result of
Theorem~\ref{thm:main-conv-result} is more general and does not
require constant stepsizes.
Another interesting observation is that the upper bounds can be
calculated without requiring any information about the
equilibrium~$\bar x$.  This is desirable, as this equilibrium is not
known beforehand and must be determined by the algorithm itself.

\section{Simulations}\label{sec:simulations}

In this section we illustrate the convergence properties of the
interconnected time-triggered system~\eqref{eq:timetrigclsys}. We
consider the IEEE 14-bus power network depicted in
Figure~\ref{fig:14-bus}, where each node has one generator and one
load according to model~\eqref{eq:swingdeltacomp}.  We assume costs at
each node $i\in[14]$ of the form
\begin{align*}
  C_i(P_{gi})=\frac12q_iP_{gi}^2+c_iP_{gi}
\end{align*}
with $q_i>0$ and $c_i\geq0$. In the original IEEE 14-bus benchmark
model, nodes $1,2,3,6,8$ have synchronous generators while the other
nodes are load nodes and have no power generation. We replicate this
by 
suitably choosing the cost at the load nodes such that the optimizer
of the economic dispatch problem \eqref{eq:ISO-OPFsc} is zero at them.
In addition, we choose $M_i\in[4,5.5]$ for generator nodes
$i\in\{1,2,3,6,8\}$ and $M_i\ll1$ for the load nodes. We set $
A_i\in[1.5,2.5],
V_i\in[1,1.06],\gain_{bi}\in[0.0005,0.001],\gain_{gi}=13.5$ for all
$i\in[14]$ and $\rho=900$. The other parameter values for the ISO
dynamics~\eqref{eq:ISOprimaldual} are $\gainl_\l=0.0004$, $\rho=3$,
$\sigma=17$. 

\tikzstyle{edge} = [draw,line width=1pt,-]
\tikzstyle{edge} = [draw,line width=1pt,-]
\tikzstyle{vertex}=[circle,
draw=black!100,line width=1pt,fill=blue!20,minimum size=16pt,inner sep=0pt]
\tikzstyle{block}=[rectangle,
draw=black!100,line width=1pt,fill=blue!30,minimum width=18mm,rounded corners]
\begin{figure}[htb]
  \centering
  \begin{tikzpicture}[scale=1.0, auto,swap]
    \definecolor{v1}{rgb}{0 ,   0.4470,    0.7410} 
    \definecolor{v2}{rgb}{    0.8500 ,   0.3250,    0.0980}
    \definecolor{v3}{rgb}{    0.9290 ,   0.6940,    0.1250}
    \definecolor{v4}{rgb}{    0.4940 ,   0.1840,    0.5560}
    \definecolor{v5}{rgb}{    0.4660 ,   0.6740,    0.1880}
    \definecolor{v6}{rgb}{    0.3010 ,   0.7450,    0.9330}
    \definecolor{v9}{rgb}{    0.6350 ,   0.0780,    0.1840}
    \definecolor{v8}{rgb}{    0.4940 ,   0.1840,    0.5560}
    \definecolor{v7}{rgb}{    1 ,   0.3,    0.3}
    \foreach \x/\y/\i/\c in {0/0/1/1,3/0/2/2,6/0/3/3,4.5/1/4/7,1.5/1/5/7,0/2/6/6,6/2/7/7,6/1/8/8,4.5/2
      /9/7,3/2/10/7,1.5/2/11/7,0/3/12/7,1.5/3/13/7,4.5/3/14/7},
    \node[vertex,fill=v\c!60] (\i) at (\x,\y) {$\i$};

    \foreach \source/\dest in {1/2,1/5,2/3,2/4,2/5,3/4,4/5,4/7,4/9,5/6,6/11,6/12,6/13,7/8,7/9,9/10,9/14,10/11,12/13,13/14},
    \path[edge] (\source) -- (\dest);
  \end{tikzpicture}
    \caption{Schematic of the modified IEEE 14-bus benchmark. Each edge
    represents a transmission line.  Red nodes represent
    loads.  All the other nodes represent synchronous generators, with
    different colors that match the ones used in
    Figures \ref{fig:hyb-sys-stab} and \ref{fig:hyb-sys-marg-stab}. The physical dynamics are modeled
    by~\eqref{eq:swingdeltacomp}.  }
  \label{fig:14-bus}
\end{figure}
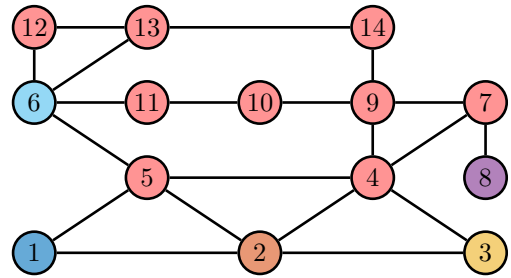

\begin{figure*}[tbh]
  \centering
  \begin{subfigure}[t]{0.32\textwidth}
    \includegraphics[width=\textwidth]{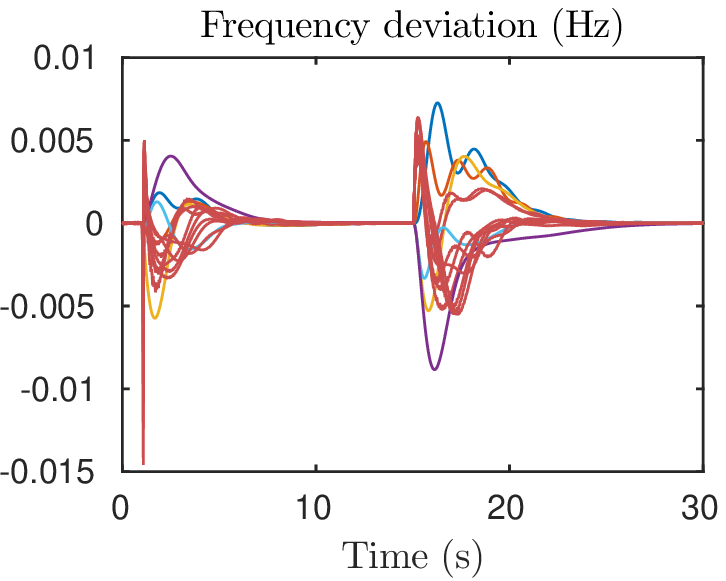}
    \caption{Evolution of the frequency deviations. After each sudden
      supply-demand imbalance, frequency is restored to its nominal
      value.}
    \label{fig:freq-s}
  \end{subfigure}
  ~ 
  \begin{subfigure}[t]{0.32\textwidth}
    \includegraphics[width=\textwidth]{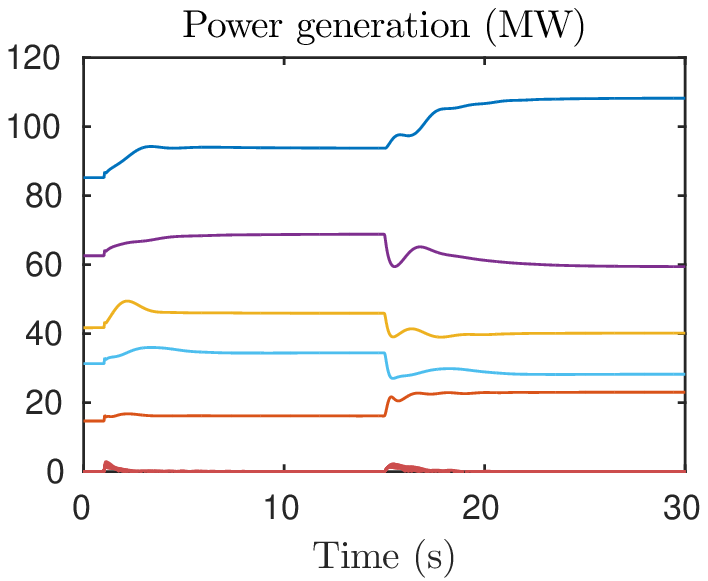}
    \caption{Evolution of the nodal power generations. After each change
      in the network, the power generation quantities converge to the
      optimal values determined by \eqref{eq:ISO-OPFsc}. }
    \label{fig:Pg-s}
  \end{subfigure}
  ~ 
  \begin{subfigure}[t]{0.32\textwidth}
    \includegraphics[width=\textwidth]{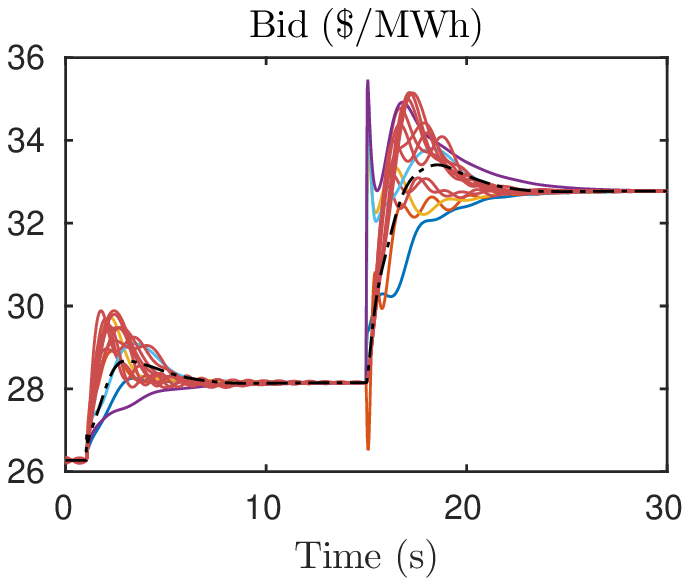}
    \caption{Evolution of the bids and the Lagrange multiplier (dashed
      black line).  As shown, the bids converge to the unique
      efficient Nash equilibrium.}
    \label{fig:bids-s}
  \end{subfigure}
  \caption{Simulations of the interconnection between the iterative
    bidding mechanism and the power network
    dynamics modeled by the time-triggered system
    \eqref{eq:timetrigclsys}. The colors in the graph corresponds to the nodes as
     depicted in Figure \ref{fig:14-bus}.  We choose identical
    inter-event times given by $t_k^l-t_{k-1}^l=\SI{2}{\milli\second},
    t_0^l-t_0^{l-1}=\SI{50}{\milli\second}$ for all $l\in\mathbb
    Z_{\geq1},k\in[25]$. As expected, the
    time-triggered system is asymptotically stable for
    sufficiently fast updates. }\label{fig:hyb-sys-stab}
\end{figure*}

\begin{figure*}[tbh]
  \centering
  \begin{subfigure}[t]{0.32\textwidth}
    \includegraphics[width=\textwidth]{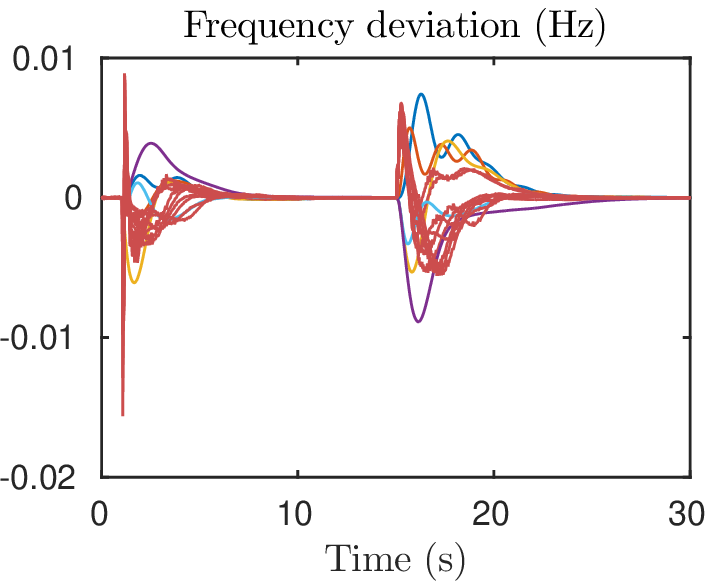}
    \caption{Evolution of the frequency deviations.}
    \label{fig:freq-rand}
  \end{subfigure}
  ~ 
  \begin{subfigure}[t]{0.32\textwidth}
    \includegraphics[width=\textwidth]{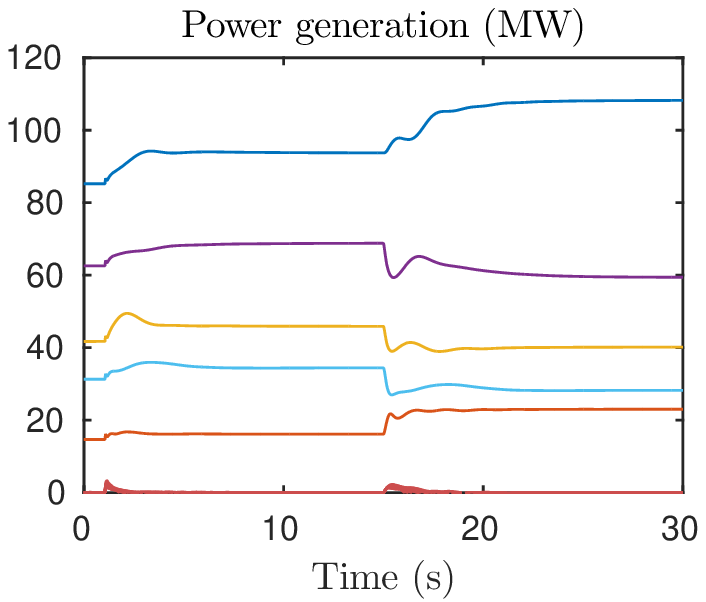}
    \caption{Evolution of each power generation.}
    \label{fig:Pg-rand}
  \end{subfigure}
  ~ 
  \begin{subfigure}[t]{0.32\textwidth}
    \includegraphics[width=\textwidth]{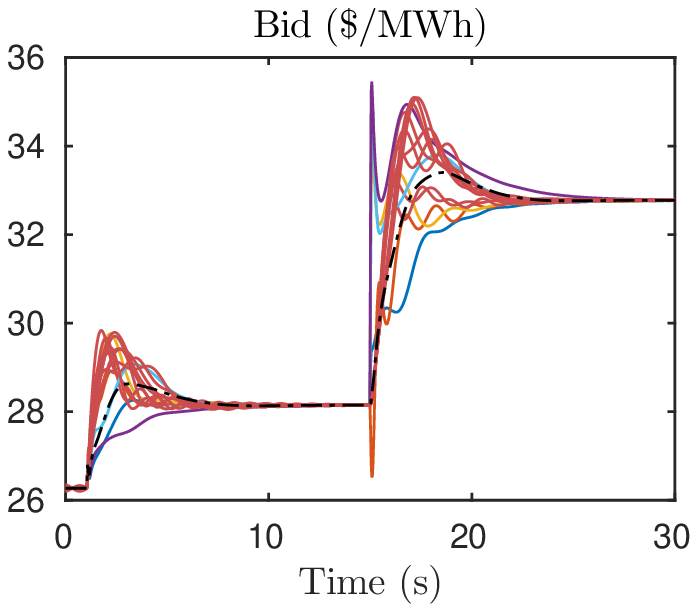}
    \caption{Evolution of the bids \& Lagrange multiplier.}
    \label{fig:bids-rand}
  \end{subfigure}
  \caption{Simulations of the time-triggered system
    \eqref{eq:timetrigclsys} under time-varying step sizes. We choose
    the time between two consecutive bid iterations randomly between
    $\SI{0.5}{\milli\second}\leq
    t_k^l-t_{k-1}^l\leq\SI{2}{\milli\second},$ for all $l\in\mathbb
    Z_{\geq1},k\in [N_l]$, and we choose the number of bid iterations
    $N_l\in\mathbb Z$ before market clearing occurs randomly in the
    interval $[20,80]$. Since the step sizes are sufficiently small, and therefore the mismatch of the time-triggered system with its continuous-time variant, the performance is similar compared to Figure \ref{fig:hyb-sys-stab}.  }\label{fig:hyb-sys-rand}
\end{figure*}

\begin{figure*}[tbh]
  \centering
  \begin{subfigure}[t]{0.32\textwidth}
    \includegraphics[width=\textwidth]{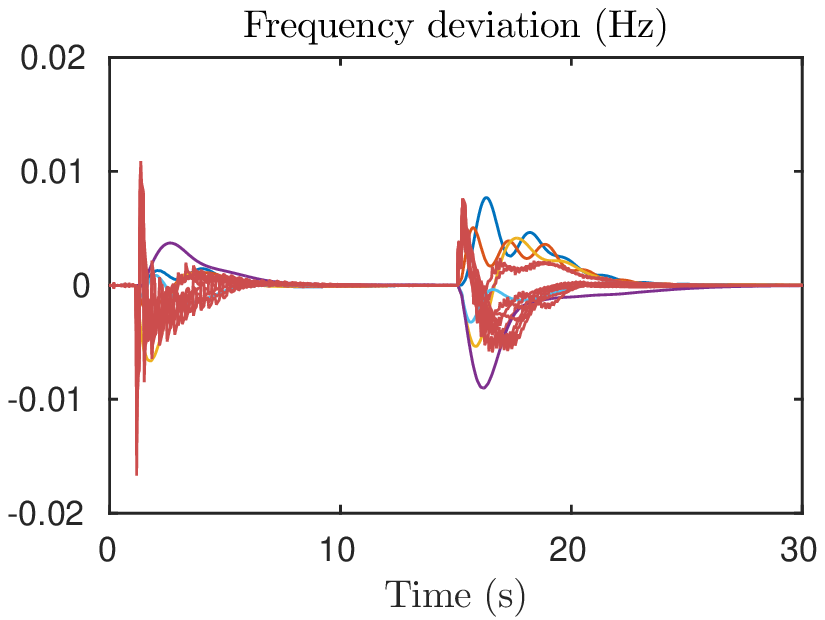}
    \caption{Compared to Figure \ref{fig:freq-s},
      there are more oscillations and a larger overshoot of the
      frequency deviations. }
    \label{fig:freq-ms}
  \end{subfigure}
  ~ 
  \begin{subfigure}[t]{0.32\textwidth}
    \includegraphics[width=\textwidth]{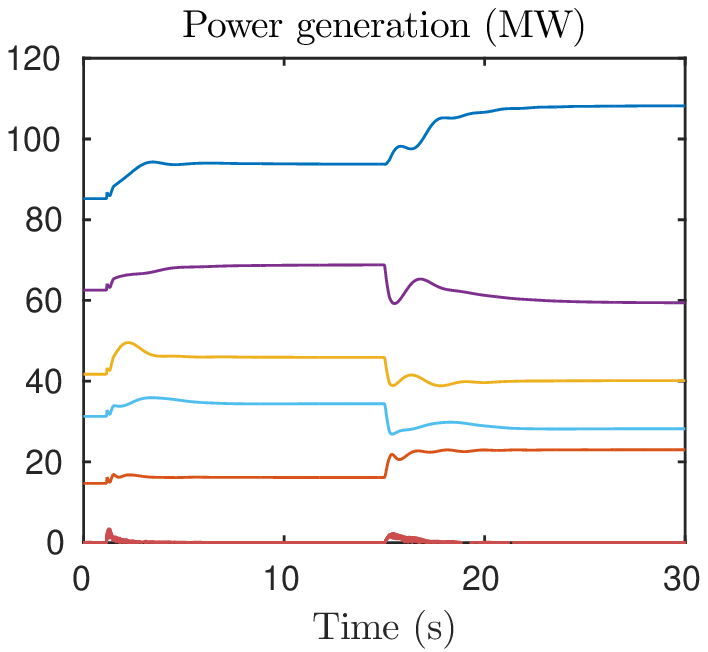}
    \caption{Evolution of the power generations at each node. }
    \label{fig:Pg-ms}
  \end{subfigure}
  ~ 
  \begin{subfigure}[t]{0.32\textwidth}
    \includegraphics[width=\textwidth]{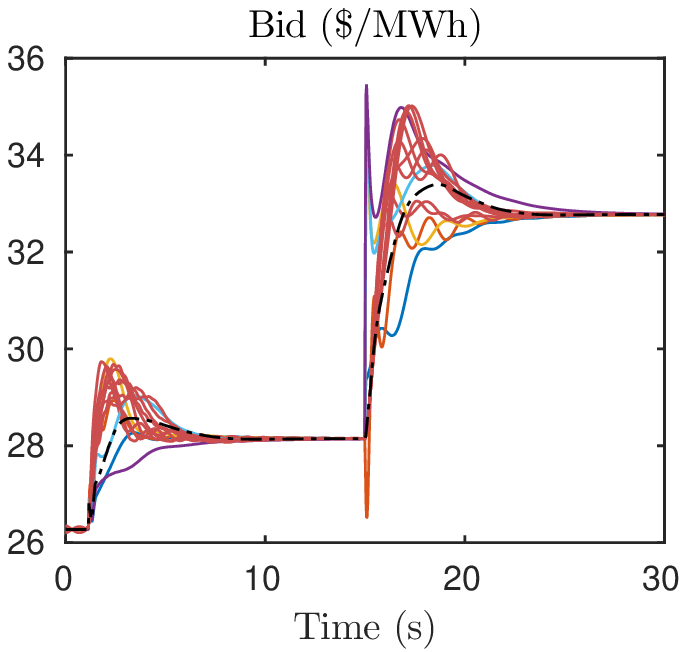}
    \caption{Evolution of the bids and the Lagrange multiplier.
      Compared to Figure  \ref{fig:bids-s}, more oscillations in the bids occur. }
    \label{fig:bids-ms}
  \end{subfigure}
  \caption{Simulations of the time-triggered
    system~\eqref{eq:timetrigclsys}. Here we consider the case
    $t_k^l-t_{k-1}^l=\SI{2}{\milli\second},
    t_0^l-t_0^{l-1}=\SI{160}{\milli\second}$ for all $l\in\mathbb
    Z_{\geq1},k\in[80]$. The scenario is the same as in Figure
     \ref{fig:hyb-sys-stab}. In this case however, the
    interconnected time-triggered system is only marginally stable; a
    small increase in either of the inter-event times results in an
    unstable system.}\label{fig:hyb-sys-marg-stab}
\end{figure*}

At time $t=\SI{0}{\second}$, the load (in \SI{}{\mega\watt}'s) is
given by
\begin{align*}
  P_d= (0, 20, 86, 43, 7, 10, 0, 0, 27, 8, 3, 6, 12, 14).
\end{align*}
Initially, we set $(q_1,q_2,q_3,q_6,q_8)=(22,128,45,60,30)$,
$(c_1,c_2,c_3,c_6,c_8)=(7.5,7.5,7.5,7.5,7.5)$ and $q_i=1500,c_i=26$
for the remaining nodes. The time-triggered system
\eqref{eq:timetrigclsys} is initialized at steady state at the
optimal generation level
\begin{align*}
  (P_{g1},P_{g2},P_{g3},P_{g6},P_{g8})=(85,15,42,31,63)
\end{align*}
and with $P_{gi}=0$ for all other
nodes. Figures~\ref{fig:hyb-sys-stab}-\ref{fig:hyb-sys-marg-stab}
depict the simulation of the time-triggered system for different
triggering times.  At $t=\SI{1}{\second}$ all the loads are increased
by 10\% and we set $c_i=28$ for the load nodes. As observed in all
figures, the trajectories converge to a new efficient equilibrium with
optimal power generation level
\begin{align*}
  (P_{g1},P_{g2},P_{g3},P_{g6},P_{g8})=(94,16,46,34,69)
\end{align*}
and $P_{gi}=0$ for all other nodes. Furthermore, at steady state the
generators all bid equal to the Lagrange multiplier which, by
Proposition~\ref{prop:exis-NE}, corresponds to an efficient Nash
equilibrium.

\begin{figure*}[tbh]
  \centering
  \begin{subfigure}[t]{0.32\textwidth}
    \includegraphics[width=\textwidth]{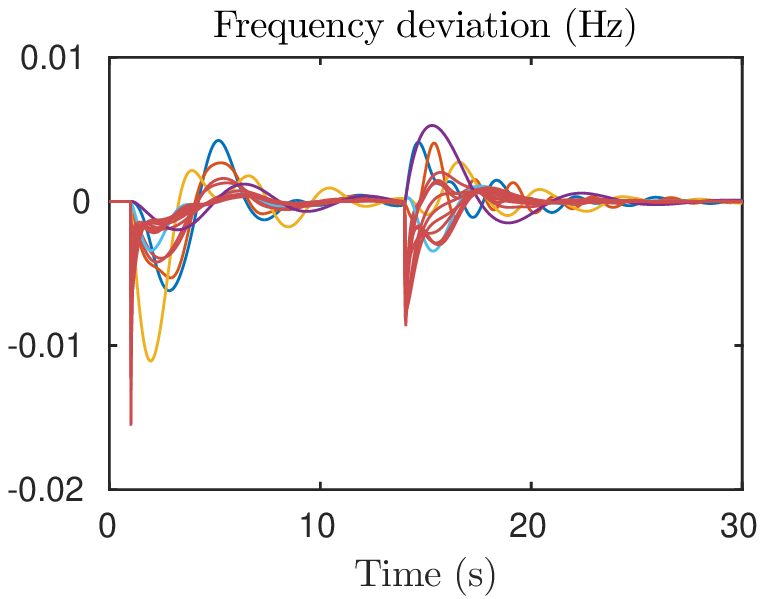}
    \caption{Evolution of the frequency deviations. Compared to Figures \ref{fig:freq-s}-\ref{fig:freq-rand}, there are more oscillations in the frequency deviations.}
    \label{fig:freq0}
  \end{subfigure}
  ~ 
  \begin{subfigure}[t]{0.32\textwidth}
    \includegraphics[width=\textwidth]{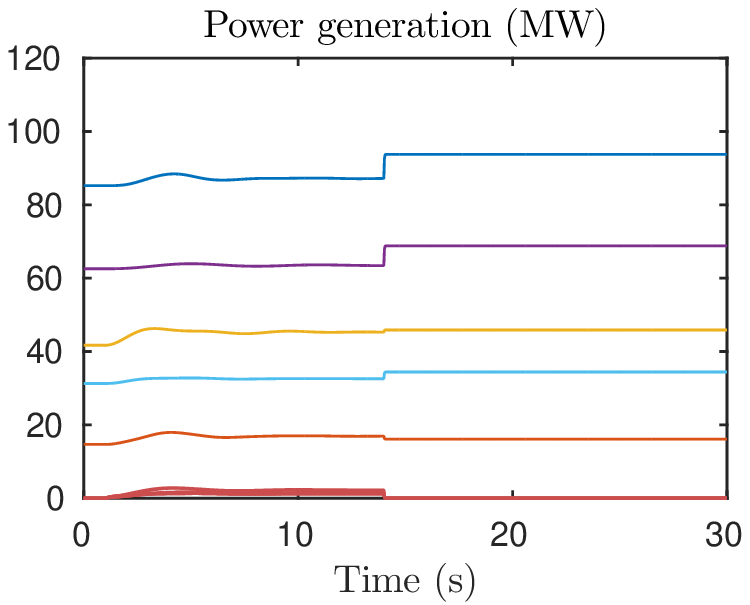}
    \caption{Evolution of each power generation. After primary and secondary controllers are activated at $t=\SI{1}{\second}$, optimal power sharing is lost.}
    \label{fig:Pg0}
  \end{subfigure}
  ~ 
  \begin{subfigure}[t]{0.32\textwidth}
    \includegraphics[width=\textwidth]{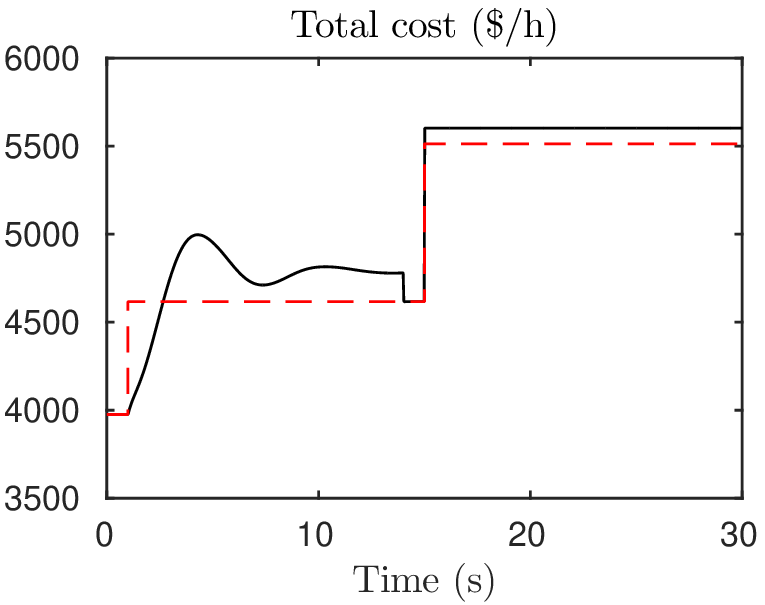}
    \caption{Evolution of the total generation costs (in black) compared to the optimal values calculated by \eqref{eq:ISO-OPFsc}. Activation of primary/secondary control, and changes in the cost function result in a loss of efficiency. }
    \label{fig:bids0}
  \end{subfigure}
  \caption{Simulations of swing equations with the primary/secondary and tertiary control layers separated. At time $t=\SI{1}{\second}$, the load is increased as in Figure \ref{fig:hyb-sys-stab} and decentralized primary/secondary controllers are activated to regulate the frequency but, as a result, optimal power sharing is lost.  At $t=\SI{14}{\second}$ the tertiary control layer is activated by resetting the setpoints optimally. After the change of the cost functions at $t=\SI{15}{\second}$, economic optimality is temporary lost again until the next time the tertiary control layer is activated (typically in the order of minutes). }\label{fig:hyb-sys-rand-sigma0}
\end{figure*}

At $t=\SI{15}{\second}$ the cost functions of the generators are
changed to $(q_1,q_2,q_3,q_6,q_8)=(23,116,48,63,38)$,
$(c_1,c_2,c_3,c_6,c_8)=(7.5,6,13.5,15,10.5)$ and $q_i=1500,c_i=33$ for
the remaining nodes. As a result, the optimal dispatch of power
changes. Due to the changes of the power generation, a temporary
frequency imbalance occurs.  As illustrated in
Figures~\ref{fig:hyb-sys-stab}-\ref{fig:hyb-sys-marg-stab}, the power
generations converge to the new optimal steady state given by
\begin{align*}
  (P_{g1},P_{g2},P_{g3},P_{g6},P_{g8})=(108, 23, 40, 28, 60).
\end{align*}
In addition, we observe that after each change of either the load or the cost function, the frequency is stabilized and the bids converge to a new efficient Nash equilibrium. The fact that the frequency transients are better in Figures \ref{fig:hyb-sys-stab}-\ref{fig:hyb-sys-rand} (with inter-event times of maximal $\SI{2}{\milli\second}$ for bidding and on average respectively $\SI{50}{\milli\second},\SI{62.5}{\milli\second}$ for market clearing) than in Figure~\ref{fig:hyb-sys-marg-stab} (with inter-event times of $\SI{2}{\milli\second}$ for bidding and $\SI{160}{\milli\second}$ for market clearing) is to be expected given the longer inter-event times in the second case.  A slight increase in the inter-event times for Figure~\ref{fig:hyb-sys-marg-stab} in either bid updating or market clearing time result in an unstable system. Figure \ref{fig:hyb-sys-rand-sigma0} illustrates the evolution of the interconnected system with the primary/secondary and tertiary control layers separated and its loss of efficiency compared to the proposed integrated design.


\section{Conclusions}\label{sec:conclusions}

This paper has studied the joint operation of the economic dispatch
and frequency regulation layers, which are traditionally separated in
the control of power networks.  The starting point of our design was a
continuous-time bid update scheme coupled with the frequency dynamics
whose equilibrium corresponds to an efficient Nash equilibrium and
zero frequency deviation. Building on the identification of a novel
LISS-Lyapunov function for this dynamics, we have characterized its
robustness properties against additive disturbances. We have exploited
the LISS-property to propose a provably correct multi-rate hybrid
implementation that combines the iterative nature of the fast bid
updates and the slower power setpoint updates with the continuous
frequency network dynamics.  Our results show that real-time iterative
bidding can successfully be interconnected with frequency control to
increase efficiency while retaining stability of the power system.
%

%
Future work will incorporate elastic demand, generator bounds, and
power flow constraints in the formulation, develop distributed and
opportunistic self-triggered implementations of the proposed dynamics,
and characterize the convergence properties of data-driven
optimization algorithms.


%
\appendices


\section{Proof of
  Theorem~\ref{thm:VeLyap}}\label{sec:proof-theor-refthm:v-1}

We structure the proof of Theorem~\ref{thm:VeLyap} in two separate
parts, corresponding to the inequalities~\eqref{eq:quadWe} and~\eqref{eq:xiein},
respectively.

\subsection{Positive definiteness of Lyapunov function $W_\e$}\label{sec:posit-defin-w_e}
Let $\bar x$ be the equilibrium of \eqref{eq:compactclsysdyn}
satisfying the hypothesis.  We now prove the existence of constants
$c_1,c_2,\e_0>0$ such that \eqref{eq:quadWe} holds, given the
constants $\e_1,\e_2,\e_3>0$. The Hessian of $W_\e$
(eq. \eqref{eq:def-We}) is given by a block-diagonal matrix
$\n^2 W_\e(x)=\blockdiag(H_1(\vp),H_2)$ with the upper left block
given by
\begin{align*}
  H_1(\vp)=
  \begin{bmatrix}
    \n^2U(\vp)&\e_0\e_1D_t^\dagger M\\
    \e_0\e_1MD_t^{\dagger T}&M
  \end{bmatrix}
\end{align*}
and the lower right block is given by
\begin{align*}
  H_2=\frac1{\sigma^2}
  \begin{bmatrix}
    \gain_b          & -\e_0\e_2\gain_g  & 0               \\
    -\e_0\e_2\gain_g & \gain_g           & -\e_0\e_3\gain_g\1 \\
    0 & -\e_0\e_3\1^T\gain_g & \gainl_\lambda
  \end{bmatrix}.
\end{align*}
We will now show that there exists sufficiently small $\e_0$ such that
$H_1(\vp),H_2$ are both positive definite for all $x\in\Omega$. To
this end, let us define the function
\begin{align}\label{eq:mathscrU}
  \mathscr U(\eta)=D_t^\dagger D\Gamma \cos(\eta)D^TD_t^{\dagger T}
\end{align}
and note that $\mathscr U(D^TD_t^{\dagger T}\vp)=\n^2U(\vp)$, implying
that $0<\mathscr U(\gamma \1)\le\n^2U(\vp)\leq \n^2U(0)=\mathscr U(0)$
for all $x\in\Omega$, see \eqref{eq:def-Omega}. Consequently, for $\D:=\e_0\e_1D_t^\dagger M$, we have
\begin{align*}
  \underbrace{
  \begin{bmatrix}
      \mathscr U(\gamma \1)&\D\\
      \D^T&M
    \end{bmatrix}}_{K_1} \leq H_1(\vp)\leq
  \underbrace{ \begin{bmatrix}
      \mathscr U(0)&\D\\
       \D^T&M
    \end{bmatrix}
         }_{K_2}, \qquad \forall x\in\Omega.
\end{align*}
%
%
By considering the Schur complements, the
matrices $K_1,H_2$ are shown to be positive definite by choosing
$\e_0>0$ sufficiently small such that
\begin{align}
  \begin{aligned}
    \mathscr U(\gamma \1)- \e_0^2\e_1^2 D_t^\dagger M D_t^{\dagger T} &>0, \\
    \gain_b-\e_0^2\e_2^2\gain_g&>0, \\
 \gainl_\lambda-\e_0^2\e_3^2\1^T\gain_b\gain_g(\gain_b-\e_0^2\e_2^2\gain_g)^{-1}\1&>0.
  \end{aligned}\label{eq:e0posdefVe}
\end{align}
Next we define
\begin{align}
  c_1&:=\min\{\lambda_{\min}(K_1),\lambda_{\min}(H_2)\},\label{eq:c1}\\
  c_2&:=\max\{\lambda_{\max}(K_2),\lambda_{\max}(H_2)\},\label{eq:c2}
\end{align}
where $\l_{\min}(A),\l_{\max}(A)$ denote the smallest and largest
eigenvalue of the matrix $A\in\mathbb R^{n\times n}$. Note that
$c_1,c_2>0$ and the following holds
\begin{align}
  0&< c_1 I\leq  \n^2 W_\e(x)\leq
     c_2 I,\qquad \forall x\in \Omega\label{eq:n2Wec1c2}
\end{align}
Note that since $W_\e(\bar x)=0, \n W_\e(\bar x)=0$, we have
\begin{align*}
  &W_\e(x)=W_\e(x)-W_\e(\bar x)\\
  &=(x-\bar x)^T\int_0^1\Big(\n W_\e((x-\bar x) \gain+\bar x)-\n W_\e(\bar x)\Big)d\gain\\
  &=(x-\bar x)^T\int_0^1\int_0^1\gain\n^2W_\e((x-\bar x)\gain\theta+\bar x)d\gain d\theta (x-\bar x).
\end{align*}
Since $\Omega$ is convex, it follows that
$x \gain \theta +(1-\gain\theta)\bar x\in \Omega$ for all
$\gain,\theta\in[0,1], x\in\Omega$. Consequently, by
\eqref{eq:n2Wec1c2} we have
\begin{align*}
  c_1I\le  \n^2 W_\e(x \gain \theta +(1-\gain\theta)\bar x)\le c_2I, \ \forall \gain,\theta\in[0,1],
\end{align*}
and $\forall x\in\Omega$.  Since $\int_0^1\int_0^1\gain d\theta
d\gain=\frac12$, inequality~\eqref{eq:quadWe}  follows.

\subsection{Dissipation inequality}\label{sec:proof-theor-refthm:v}

Here we establish the inequality~\eqref{eq:xiein}. First we consider
the case without disturbance, i.e., $d=0$. Given the equilibrium
$\bar x$ of \eqref{eq:compactclsysdyn}, we define $\tilde x:=x-\bar x$
and likewise $\tilde \vp,\tilde \w,\tilde b,\tilde P_g,\tilde
\lambda$.  Then, the system~\eqref{eq:compactclsysdyn} reads
as 
\begin{align*}
  \dot{\tilde\vp}&=D_t^T\tilde \w,
  \\
  M\dot{\tilde \w}&=-D_t(\nabla U(\vp)-\nabla U(\bar\vp))-A\tilde \w+\tilde P_g,\\
  \gain_b \dot{\tilde b}&=\tilde P_g-Q^{-1}\tilde b,\\
  \gain_g\dot{\tilde P}_g&=\1\tilde \l-\tilde b-\rho \1\1^T\tilde P_g -\sigma^2 \tilde \w,\\
  \gainl_\l\dot{\tilde \l}&=-\1^T\tilde P_g.
\end{align*}
In addition, note that $W_\e$ (eq. \eqref{eq:def-We}) takes the form
\begin{align}
  &W_\e(x)=V(x)+V_\e(x),\\
  &V_\e(x)=\e_0\e_1\tilde \vp^TD_t^{\dagger}M\w-\frac{\e_0\e_2}{\sigma^2}
    \tilde b^T\gain_g\tilde P_g- \frac{\e_0\e_3}{\sigma^2}
    \tilde \l\1^T\gain_g\tilde P_g.\label{eq:Ve-recall}
\end{align}
Next, we determine the time-derivative of the individual terms of the
candidate Lyapunov function $W_\e$.

\noindent \textbf{(0):} 
First, observe from \eqref{eq:V-diss-ineq} that
\begin{align*}
  \dot V&= - \w^TA\w-\frac{1}{\sigma^2}(b-\bar b)^TQ^{-1}(b-\bar b)\\&-\frac\rho{\sigma^2}(P_g-\bar P_g)^T\1\1^T(P_g-\bar P_g).
\end{align*}
%
%

\noindent \textbf{(1):} 
The time-derivative of the first term of $V_\e$ satisfies
\begin{align*}
  &\ddt \tilde\vp^TD_t^{\dagger}M\w=\tilde\w^TMD_t^{\dagger T}D_t^T\tilde \w\\
  &-\tilde \vp^T D_t^{\dagger}D_t(\nabla U(\vp)-\nabla U(\bar\vp))-\tilde \vp^TD_t^\dagger A\tilde \w+\tilde \vp^TD_t^\dagger \tilde P_g.
\end{align*}
By exploiting $D_t^{\dagger}D_t=I$, the second term is rewritten as
\begin{align*}
  -\tilde \vp^T D_t^{\dagger}D_t(\nabla U(\vp)-\nabla U(\bar\vp))
  &=-\tilde \vp^T\mathcal U(\vp)\tilde\vp^T
\end{align*}
where we used that $\n U(\vp)-\n U(\bar\vp)=\mathcal U(\vp)(\vp-\bar \vp)$ with
\begin{align}\label{eq:mathcalU}
  \mathcal U(\vp)=\int_0^1\n^2 U((\vp-\bar \vp) \theta+\bar\vp)d\theta.
\end{align}
Since
$\mathcal U(\vp)\geq \mathscr U(\1 \gamma)=D_t^\dagger D\Gamma
\cos(\1\gamma)D^TD_t^{\dagger T}$ (see eq. \eqref{eq:mathscrU}) for
all $x\in\Omega$, we obtain
\begin{align*}
  \ddt \tilde\vp^TD_t^{\dagger}M\w&\leq\tilde\w^TMD_t^{\dagger T}D_t^T\tilde \w-\tilde \vp^T\mathscr U(\1\gamma)\tilde\vp^T\\&-\tilde \vp^TD_t^\dagger A\tilde \w+\tilde \vp^TD_t^\dagger \tilde P_g.
\end{align*}

\noindent \textbf{(2):} 
For the second term of $V_\e$ the following holds:
\begin{align*}
  \ddt \tilde b^T\gain_g\tilde P_g&= \tilde P_g^T \gain_{gb} \tilde P_g-\tilde P_g^T \gain_{gb}Q^{-1}\tilde b+\tilde b^T\1\tilde \l\\&-\tilde b^T\tilde b-\rho \tilde b^T\1\1^T\tilde P_g -\sigma^2 \tilde b^T\tilde \w,
\end{align*}
where we define $\gain_{gb}:=\gain_g\gain_b^{-1}$.

\noindent \textbf{(3):} 
Similarly, by defining $\gain_{g\l}:=\gain_g\gain_{\l}^{-1}$ we obtain
\begin{align*}
  \ddt \tilde \l\1^T\gain_g\tilde P_g&=-\tilde P_g^T\gain_{g\l} \1\1^T\tilde P_g+ n\tilde\l^2-\tilde\l \1^T\tilde b\\&-\rho n \tilde\l \1^T\tilde P_g -\sigma^2\tilde\l \1^T \tilde \w.
\end{align*}
By combining the above calculations, we can show that the time-derivative of $W_\e$ satisfies
\begin{align*}
  \dot W_\e&=\dot V+\dot V_\e\leq \frac12\e_0(x-\bar x)^T\mathcal P^T\Xi_\e\mathcal P(x-\bar x).
\end{align*}
where $\Xi_\e$ is given by \eqref{eq:Theta} (see page
\pageref{eq:Theta}), $\mathcal P$ takes the form
\begin{align*}
  \mathcal P=\begin{bmatrix}
    0 & I & 0              & 0              & 0               \\
    0 & 0 & \frac1\sigma I & 0              & 0               \\
    0 & 0 & 0              & \frac1\sigma I & 0               \\
    0 & 0 & 0              & 0              & \frac1\sigma   \\
    I & 0 & 0              & 0              & 0                    
  \end{bmatrix},
\end{align*}
and
$\mathscr M:=MD_t^{\dagger T}D_t^T+D_tD_t^{\dagger}M, \mathscr
T:=\gain_{g\l}\1\1^T+\1\1^T\gain_{g\l}$.  Next, we will show that there
exists $\e_0,\e_1,\e_2,\e_3>0$ such that $\Xi_\e$ is positive
definite.
\begin{figure*}
  \centering
  \begin{align}\label{eq:Theta}
  \Xi_\e&=
  \bbordermatrix{
                 & \w                                                      & b/\sigma                      & P_g/\sigma                                                    & \l/\sigma       & \vp                          \cr
\w               & \frac{2}{\e_0}A-\e_1\mathscr M & -\e_2\sigma I                 & 0                                                             & -\e_3 \sigma \1 & \e_1AD_t^{\dagger T}         \cr
\frac b\sigma    & -\e_2\sigma I                                           & -2 \e_2I+\frac{2}{\e_0}Q^{-1}       & -\e_2(Q^{-1}\gain_{gb}+\rho \1\1^T)                                 & (\e_2-\e_3)\1   & 0                            \cr
\frac{P_g}\sigma & 0                                                       & -\e_2(\gain_{gb}Q^{-1}+\rho \1\1^T) & 2\e_2\gain_{gb}+\frac{2}{\e_0}\rho \1\1^T-\e_3\mathscr T & -\e_3 n\rho \1  & -\e_1\sigma D_t^{\dagger T}  \cr
\frac\l\sigma    & -\e_3 \sigma \1^T                                       & (\e_2-\e_3)\1^T               & -\e_3n \rho \1^T                                              & 2n\e_3          & 0                            \cr
\vp              & \e_1D_t^{\dagger}A                                      & 0                             & -\e_1\sigma D_t^{\dagger}                                     & 0               & 2 \e_1 \mathscr U(\1 \gamma) \cr
  }
\end{align}
\end{figure*}
This can be done by successive use of the Schur complement. 
In particular, for $A\in\mathbb R^{n\times n},B\in\real^{n\times m},C\in\real^{m\times m},\beta>0$, recall that
\begin{align*}
  \begin{bmatrix}
     \beta A&B\\
    B^T&C
  \end{bmatrix}>0\quad \Longleftrightarrow\quad C>0\ \ \&\ \   \beta A-BC^{-1}B^T>0.
\end{align*}
For successively applying this result to $\Xi_\e$, given by
\eqref{eq:Theta}, let us first fix $\e_1,\e_3>0$. Then $\e_2$ can be
chosen sufficiently large such that lower-right $3\times3$ block
submatrix of $\Xi_\e$ is positive definite. Then we can choose a
$\e_0>0$ sufficiently small such that \eqref{eq:e0posdefVe} holds and
$\Xi_\e>0$. Here, note that choosing $\e_0$ smaller does not affect
the positive definiteness of the lower-right $3\times 3$ block
submatrix of $\Xi_\e$.  By construction of $\e_0,\e_1,\e_2,\e_3$,
there exist constants $c_1,c_2\in\mathbb R_{>0}$ such that
\eqref{eq:quadWe} holds for all $x\in\Omega$, see also Section
\ref{sec:posit-defin-w_e}. In addition, for this choice of $\e$ we
have that $\Xi_\e>0$ and, as a result, there exists
$\hat \alpha:=\frac12\e_0\l_{\min}(\mathcal P^T\Xi_\e\mathcal P)>0$ such
that 
\begin{align*}
  (\nabla W_\e(x))^TF(x)&\leq -\hat \alpha\|x-\bar
  x\|^2
\end{align*}
for all $x\in\Omega$.  Next, we consider the case when the disturbance
is present. Let $\chi$ satisfy $0<\chi<\hat \alpha/(L_W\norm
B)$. Then, by exploiting the Lipschitz property of $\n W_\e$,
\begin{align*}
  &(\nabla W_\e(x))^T(F(x)+Bd)\leq -\hat \alpha\|x-\bar
  x\|^2+\nabla W_\e(x))^TBd\\
  &\leq -\hat \alpha\|x-\bar x\|^2+L_W\norm B\norm{x-\bar x}\norm
  d\\& \leq -(\hat \alpha-L_W\norm B\chi)\norm{x-\bar
      x}^2=-\alpha\norm{x-\bar x}^2
  \end{align*}
  with $\alpha:=\hat \alpha-L_W\norm B\chi)>0$ and thus
  \eqref{eq:xiein} holds.  This concludes the proof of Theorem
  \ref{thm:VeLyap}. \hfill $\blacksquare$

\section{Proof of
  Theorem~\ref{thm:main-conv-result}}\label{sec:two-time-triggering}

Here we prove Theorem~\ref{thm:main-conv-result}. To do so, we rely on
\gwin, which in general allows to bound the evolution of
continuous-time and discrete-time signals described by differential
and difference equations, respectively. Given the hybrid nature of the
time-triggered dynamics~\eqref{eq:timetrigclsys}, we rely on a version
of \gwin\ for hybrid systems developed in~\cite{noroozi2014gronwall}.
Adapted for our purposes, it states the following.
\begin{proposition}\longthmtitle{Generalized \gwin\
    \cite{noroozi2014gronwall}}\label{prop:hybridGronwall}
  Let $t \mapsto y(t) \in \real$ be a continuous signal, $t \mapsto p(t)
  \in \real$ be a continuously differentiable signal,
  $r:=\{r_j\}_{j=0}^{k-1}$ be a nonnegative sequence of real
  numbers, $q\geq0$ a constant, and $E:=\{t_j\}_{j=0}^{k+1}$, $k \in
  \mathbb Z_{\ge 0}$ be a sequence of times satisfying $t_j
  < t_{j+1}$ for all $j \in \{0, \dots, k\}$. Suppose that for
  all $t\in [t_0,t_{k+1}]$, the elements $y$, $p$, and $r$ satisfy
  \begin{align*}
    y(t)\leq p(t)+q\myint_{t_0}^t y(s)ds+\mysum[]_{j=0}^{i(t)-1}r_jy(t_{j+1})
  \end{align*}
  with $i(t):=\max\{i\in\mathbb Z_{\ge0}: t_i\leq t, t_{i}\in E\}$ 
  for $t<t_{k+1}$ and $i(t_{k+1}):=k$. Then,
  \begin{align*}
    y(t)&\leq p(t_0)h(t_0,t)+\myint_{t_0}^t h(s,t)\dot p(s)ds
  \end{align*}
  for all  $t\in[t_0,t_{k+1}]$ where for all $t_0 \le s \le t \le t_{k+1}$, 
  \begin{align*}    
    h(s,t) :=\exp\big(q(t-s)+\mysum_{j=i(s)}^{i(t)-1}\log(1+r_j)\big).
  \end{align*}
\end{proposition}

We are now ready to prove Theorem~\ref{thm:main-conv-result}.

\subsubsection*{Proof of Theorem~\ref{thm:main-conv-result}}

Let $\{ \{t_k^l\}_{k=0}^{N_l} \}_{l=0}^\infty $ be a sequence of times
satisfying the hypotheses. Consider a trajectory $t \mapsto x(t)$
of~\eqref{eq:timetrigclsys} with $x(0)$ belonging to a neighborhood of
$\bar x$. The definition of this neighborhood will show up later. Our
proof strategy involves showing the monotonic decrease of the function
$W_\e$ (cf.~\eqref{eq:def-We}) along this arbitrarily chosen
trajectory.  Consider any $t \in \realnonnegative$ such that $t \not
\in \{t_k^l\}_{k=0}^{N_l}$ for any $l \in \mathbb Z_{\ge 0}$ and
$x(t)\in\Omega$ where $\Omega$ is defined by
\eqref{eq:def-Omega}. With a slight abuse of notation let $l$ and $k
\in \{0,\dots,N_l -1\}$ be fixed such that $t \in (t_k^l,
t_{k+1}^l)$. Then, using the expression of $F(x)= f(x)+g(x)+h(x)$
given in~\eqref{eq:gennonsys}, one can write the evolution of $x$ at
$t$ for the considered trajectory as
  \begin{align*}
    \dot x(t)=F(x(t))+g(x(t_k^l))-g(x(t))+h(x(t_0^l))-h(x(t)).
  \end{align*}
  \textbf{(I) Dissipation inequality:} Note that at $t$ the evolution of
  $W_\e$ is equal to the dot product between the gradient of $W_\e$ and
  right-hand side of the above equation. Hence, we get 
  \begin{align}
    \dot W_\e(x(t)) & = \n W_\e(x(t))^\top \Bigl( F(x(t))
    +g(x(t_k^l))-g(x(t)) \notag
    \\
    & \qquad +h(x(t_0^l))-h(x(t)) \Bigr). \label{eq:lie-W-eps}
  \end{align}
  From~\eqref{eq:xiein}, we have $\n W_\e(x(t))^\top F(x(t)) \le -
  \alpha \norm{x(t) - \bar x}^2$. Moreover, since maps $\n W_\e$, $g$, and
  $h$ are globally Lipschitz and $\n W_\e(\bar x) = 0$, one has
  $\norm{\n W_\e(x(t))} \le L_W \norm{x(t) - \bar x}$,
  $\norm{g(x(t_k^l)) - g(x(t))}  \le L_g \norm{ x(t_k^l) - x(t)}$, and
  $\norm{h(x(t_0^l))-h(x(t))}  \le L_h \norm{x(t_0^l) - x(t)}$.  Using
  these bounds in~\eqref{eq:lie-W-eps}, we get
  \begin{align}
    \dot W_\e(x(t)) & \leq -\alpha \norm{x(t)-\bar x}^2 +L_W
    \norm{x(t)-\bar x}\Bigl( L_g \|x(t) \notag
    \\
    & \qquad -x(t_k^l)\| +L_h \norm{x(t)-x(t_0^l)} \Bigr).
    \label{eq:dotWe-time-trig}
  \end{align}
  Next, we provide bounds on $\|x(t)-x(t_k^l)\|$ and $\|x(t)-x(t_0^l)\|$
  in terms of $\|x(t)-\bar x\|$, $t-t_k^l$, and $t-t_0^l$. 
  To reduce the notational burden, we drop the superscript $l$ from the
  time instances $\{t_i^l\}_{i=1}^{N_l}$. In addition, we define 
  \begin{alignat*}{2}
    x_k &:=x(t_k), &\quad \zeta_k(t)&:=t-t_k,
    \\
    \zeta_j^k&:=\zeta_j(t_k)=t_k-t_j, &\qquad
    \xi^l(t)&:=\zeta_0(t)=t-t_0.
  \end{alignat*}
  %
  %
  \textbf{(II) Bounds on $\mathbf{\norm{x(t)-x(t_k^l)}}$:} Note that
  $x(t)$ can be written using~\eqref{eq:gentwotrigsys} as the line
  integral
  \begin{align}
    x(t)- & x_k  =  \myint_{t_k}^tf(x(s))ds+\zeta_k(t)g(x_k)+\zeta_k(t)h(x_0)
    \notag
    \\
    & =  \myint_{t_k}^t(f(x(s))-f(x_k))ds +\zeta_k(t)(f(x_k)-f(\bar x))
    \notag
    \\
    & \quad  +\zeta_k(t)(g(x_k)-g(\bar x)+h(x_0)-h(\bar x)). \label{eq:xt-xkeq}
  \end{align}
  Above, we have added and subtracted $\zeta_k(t) f(x_k)$
  and subtracted $f(\bar x) + g( \bar x) + h( \bar x)$ as $\bar x$ is an
  equilibrium.  Using Lipschitz bounds and triangle inequality
  in~\eqref{eq:xt-xkeq} we obtain 
  \begin{align}
    \norm{x(t)&-x_k} \leq L_f \myint_{t_k}^t \norm{x(s)-x_k} ds
    \label{eq:state-evol-1}
    \\
    & +\zeta_k(t)(L_f+L_g) \norm{x_k-\bar x}+\zeta_k(t)L_h
    \norm{x_0-\bar x}.  \notag 
  \end{align}
  From above, we wish to obtain an upper bound on $ \norm{x(t)-x_k}$ 
  that is independent of the state at times $s \in
  (t_k,t)$. To this end, we employ Gronwall's inequality as stated in a
  general form in Proposition~\ref{prop:hybridGronwall}. Drawing a
  parallelism between the notations, for~\eqref{eq:state-evol-1}, we
  consider $E = \emptyset$, $r=0, y(t)=\|x(t)-x_k\|, q=L_f,
  p(t)=\zeta_k(t)(L_f+L_g)\|x_k-\bar x\|+\zeta_k(t)L_h\|x_0-\bar x\|$.
  Then, applying Proposition~\ref{prop:hybridGronwall} and integrating
  the then obtained right-hand side yields  
  \begin{align}\label{eq:x(t)-xkineq}
  \begin{aligned}
    \|x(t)-x_k\|\leq \Bigl(1+\tfrac{L_g}{L_f} \Bigr) &\|x_k-\bar
    x\|(e^{L_f\zeta_k(t)}-1)
    \\
    +\tfrac{L_h}{L_f}&\|x_0-\bar x\|(e^{L_f \zeta_k(t)}-1).
  \end{aligned}
  \end{align} 
  Bounding the above inequality using the triangle inequality 
    $\norm{x_k-\bar x} \leq \norm{x(t)-x_k} + \norm{x(t)-\bar x}$, 
  collecting coefficients of $\norm{x(t) - x_k}$
  on the left-hand side, and rearranging gives 
  \begin{align}
    & \norm{x(t) - x_k} \le \frac{L_h(e^{L_f\zeta_k(t)}-1)}
    {L_f-(L_f+L_g)(e^{L_f\zeta_k(t)}-1)}\|x_0-\bar x\| \notag
    \\
    & \quad
    +\frac{(L_f+L_g) (e^{L_f\zeta_k(t)}-1) }
    {L_f-(L_f+L_g)(e^{L_f\zeta_k(t)}-1)} \|x(t)-\bar x\|.
    \label{eq:xt-xk-bound}
  \end{align}
  \textbf{(III) Bounds on $\mathbf{\norm{x(t)-x(t_0^l)}}$:}  Our next
  step is to provide an upper bound on the term $\norm{x(t) - x_0}$.
  Recall that the considered trajectory
  satisfies~\eqref{eq:gentwotrigsys} and so, the line integral over the
  interval $[t_0,t]$ gives
  \begin{align*}
    x(t)-x_0 & = \myint_{t_0}^t f(x(s))ds + \textstyle \mysum[]_{j=0}^{k-1} \zeta_j^{j+1}
    g(x_j) 
    \\
    & \qquad + \zeta_k(t)g(x_k) + \xi^l(t)h(x_0).
  \end{align*}
  As done before, on the right-hand side, we add and subtract the
  terms $\xi^l(t) f(x_0)$ and $\xi^l(t) g(x_0)$ and then subtract
  $f(\bar x) + g(\bar x) + h(\bar x)$. This gives us
  \begin{align*}
    x(t) &-  x_0  = \myint_{t_0}^t(f(x(s))-f(x_0))ds 
    \\
    & \textstyle + \mysum[]_{j=0}^{k-1} \zeta_{j}^{j+1}(g(x_j)-g(x_0)) + \zeta_k(t)
    (g(x_k) - g(x_0)) 
    \\
    & + \xi^l(t)
    (f(x_0) - f(\bar x) + g(x_0) - g(\bar x) + h(x_0) - h(\bar x) )
  \end{align*}
  By defining $L:=L_f+L_g+L_h$, taking the norms, using the global
  Lipschitzness, we obtain from above 
  \begin{align*}
    \|x(t) & -x_0\| \leq
    L_f\myint_{t_0}^t\|x(s)-x_0\|ds + \xi^l(t)L\|x_0-\bar x\|
    \\
    & \textstyle + L_g\mysum[]_{j=0}^{k-1}\zeta_{j}^{j+1}\|x_{j}-x_0\| + L_g \zeta_k(t)
    \norm{x_k - x_0}.
  \end{align*}
  Consider any $\hat t \in [t,t_{k+1}]$ and note that $\zeta_k(t) \le
  \zeta_k(\hat t)$. Using this bound and the fact that the first term in
  the above summation is zero, we write 
  \begin{align*}
    \|x(t)- & x_0\|  \leq
    L_f\myint_{t_0}^t\|x(s)-x_0\|ds + \xi^l(t)L\|x_0-\bar x\|
    \\
    & \, \textstyle + L_g\mysum[]_{j=0}^{k-2}\zeta_{j+1}^{j+2}\|x_{j+1}-x_0\| + L_g
    \zeta_k(\hat t) \norm{x_k - x_0}.  
  \end{align*}
  We now apply Proposition~\ref{prop:hybridGronwall} to give a bound for
  the left-hand side independent of $x(s)$, $s \in
  (t_0,t]$. In order to do so, the elements corresponding to those in
  the Gronwall's inequality are: $E=\{t_j\}_{j=0}^{k}\cup  \{\hat t\}$,
  $y(t)=\|x(t)-x_0\|, p(t)=\xi^l(t)L\|x_0-\bar x\|,q=L_f,
  r_j=L_g\zeta_{j+1}^{j+2}$ for $j=0,\ldots,k-2$, and $r_{k-1}=\hat
  t-t_k$.  From Proposition~\ref{prop:hybridGronwall}, we get
  \begin{align}
    \|x(t)-x_0\| \leq  L \|x_0-\bar x\| \myint_{t_0}^{t} h(s,t) ds,
    \label{eq:gronwall-imply}
  \end{align}
  where
  $ h(s,t) = \exp\big(\textstyle\myint_s^{t}L_fd\gain + \textstyle
  \mysum[]_{j=i(s)}^{k-2}\log (1+\zeta_{j+1}^{j+2}L_g) + \log(1+L_g
  \zeta_k(\hat t)) \big)$
    and $i(s)$ is as defined in Proposition~\ref{prop:hybridGronwall}.
  Using $\log(1+x) \le x$ for $x \ge 0$ and the fact
  that the exponential is a monotonically increasing function, we get
  \begin{align*}
   \textstyle  h(s,t) \le  \exp \left( L_f (t-s) +    L_g \mysum[]_{j=i(s)}^{k-2}
    \zeta_{j+1}^{j+1} + L_g \zeta_k(\hat t) \right).
  \end{align*}
  By noting that $s \le i(s) + 1$ and $t \le \hat t$, we can upper
  bound the right-hand side as
  $h(s,t) \le \exp \left( L_f (t-s) + L_g (\hat t-s) \right)$.  Since
  $\hat t$ was chosen arbitrarily in the interval $[t,t_{k+1}]$, we
  pick it equal to $t$. Thus,
  $h(s,t) \le \exp\left( (L_g+L_f) (t-s) \right)$.  Substituting this
  inequality in~\eqref{eq:gronwall-imply} yields
  \begin{align}
    \|x(t)& -x_0\| \leq L \norm{x_0 - \bar x}
    \myint_{t_0}^{t}e^{(L_f+L_g)(t-s)}ds \notag
    \\
     & =  \tfrac{L}{L_f+L_g}(e^{(L_f+L_g)\xi^l(t)}-1)\|x_0-\bar x\|.
     \label{eq:xt-xo-bound}
  \end{align}
  This inequality when used in the right-hand side of the triangle
  inequality $\|x_0-\bar x\|\leq \|x(t)-x_0\|+\|x(t)-\bar x\|$ yields
  after rearrangement the following
  \begin{align}
    \|x_0-\bar x\|
    &\leq
      \frac{L_f+L_g}{L_f+L_g-L(e^{(L_f+L_g)\xi^l(t)}-1)}\|x(t)-\bar
      x\|.\label{eq:x0-bx-bound}
  \end{align}
  Subsequently, using the above bound in~\eqref{eq:xt-xo-bound} gives
  \begin{align}
    \|x(t)-x_0\|
    &\leq
      \frac{L(e^{(L_f+L_g)\xi^l(t)}-1)}{L_f+L_g-L(e^{(L_f+L_g)\xi^l(t)}-1)}
      \|x(t)-\bar x\|.\label{eq:xt-x0-bound}
  \end{align}
  Combining inequalities~\eqref{eq:xt-xk-bound}
  and~\eqref{eq:x0-bx-bound} we obtain
  \begin{align}
    \norm{x(t) & - x_k} \le \frac{L_h(e^{L_f\zeta_k(t)}-1)}
    {L_f-(L_f+L_g)(e^{L_f\zeta_k(t)}-1)} \cdot \notag
    \\
    & \quad \frac{L_f+L_g}{L_f+L_g-L(e^{(L_f+L_g)\xi^l(t)}-1)}\|x(t)-\bar
    x\|\notag
    \\
    & 
    +\frac{(L_f+L_g) (e^{L_f\zeta_k(t)}-1) }
    {L_f-(L_f+L_g)(e^{L_f\zeta_k(t)}-1)} \|x(t)-\bar x\|
    \label{eq:xt-xk-bound2}
  \end{align}
  \textbf{(IV) Monotonic decrease of $\mathbf{W_\e}$:} Note first that
  following~\eqref{eq:xt-x0-bound} and using the bound
  $\xi^l(t) \le \bar \xi$ yields
  \begin{align*}
    \norm{x(t) - x_0} \le
    \frac{L(e^{(L_f+L_g)\bar \xi}-1)}{L_f+L_g-L(e^{(L_f+L_g)\bar
    \xi}-1)} \norm{x(t)-\bar x}.
  \end{align*}
  Using the definition of $\bar \xi$, one gets
  $e^{(L_f + L_g) \bar \xi} - 1 = \frac{\beta( L_f + L_g)}{L (L_W L_h
    + \beta)}$. Substituting this value in the above inequality and
  simplifying the expression provides us
  \begin{align}
    \norm{x(t) - x_0} \le (\beta/(L_W L_h) )
    \norm{x(t) - \bar x}.
  \label{eq:xt-x0-bound-final}
  \end{align}
  In a similar way, using the bound on $\xi^l(t)$ and substituting the
  value of $e^{(L_f + L_g) \bar \xi} - 1$ in~\eqref{eq:xt-xk-bound2}
  and simplifying yields
    \begin{align*}
      \norm{x(t) - x_k}
      & \le \frac{(e^{L_f\zeta_k(t)}-1)(L +
        \beta/L_W)}
        {L_f-(L_f+L_g)(e^{L_f\zeta_k(t)}-1)} \norm{x(t) - \bar x}.
  \end{align*}
  Note that $\zeta_k(t) \le \bar \zeta$. Using this bound and the
  definition of $\bar \zeta$ in the above inequality gives
  \begin{align}
    \norm{x(t) - x_k} \le  \tfrac{\alpha - \beta}{L_W L_g} 
    \norm{x(t) - \bar x}.
   \label{eq:xt-xk-b-final}
  \end{align}
  Finally, substituting~\eqref{eq:xt-x0-bound-final}
  and~\eqref{eq:xt-xk-b-final} in~\eqref{eq:dotWe-time-trig} and using
  the fact that $\beta < \alpha$, we obtain $\dot W_\e (x(t)) <
  0$. Recall that $t \in \realnonnegative$ was chosen arbitrarily
  satisfying $t \not \in \{t_k^l\}_{k=1}^{N_l}$ for any
  $l \in \mathbb Z_{\ge 0}$. Therefore, $W_\e$ monotonically decreases
  at all times along the trajectory except for a countable number of
  points. Further, the map $t \mapsto W_\e(x(t))$ is continuous.
  Therefore, we conclude that the trajectory initialized in a compact
  level set of $W_\e$ contained in $\Omega$ converges asymptotically
  to the equilibrium point $\bar x$. This completes the proof.  \hfill
  $\blacksquare$

\bibliographystyle{IEEEtran}
\bibliography{IEEEabrv,bibfile}

\end{document}